[1994/12/01]
\documentclass[12pt]{article}
\title{Towards $\psi-$extension of Finite Operator Calculus of Rota}
\author{A.K.Kwa\'sniewski \\  
\\ Institute of Computer Science, Bia{\l}ystok University\\
PL-15-887 Bia{\l}ystok, ul.Sosnowa 64, POLAND\\
Higher School of Mathematics and Applied Computer Science,\\
PL-15-364 Bia{\l}ystok, Czysta 11\\
\underline {e-mail: kwandr@noc.uwb.edu.pl}}

\usepackage{amsmath,amsthm}

\chardef\bslash=`\\ 
\newtheorem{thm}{Theorem}[section]
\newtheorem{cor}{Corollary}[section]

\newtheorem{prop}{Proposition}[section]
\newtheorem{com}{Comment}[section]
\newtheorem{ex}{Example}[section]
\newtheorem{rem}{Remark}[section]
\newtheorem{obs}{Observation}[section]
\newtheorem{conc}{Conclusion}[section]
\newtheorem{defn}{Definition}[section]

\begin{document}
\maketitle

\begin{abstract}

$\psi $- extension of Gian-Carlo Rota`s finite operator calculus due to
Viskov \cite{1,2} is further developed. The extension relies on the
notion of $\partial _{\psi}  $-shift invariance and $\partial _{\psi}
$-delta operators. Main statements of Rota`s finite operator calculus are
given their $\psi $-counterparts. This includes Sheffer $\psi $-polynomials
properties and Rodrigues formula - among others. Such $\psi $-extended
calculus delivers an elementary umbral underpinning for $q$-deformed quantum
oscillator model and its possible generalisations.

$\partial _q$-delta operators and their duals and similarly $
\partial _\psi $-delta operators with their duals are pairs of generators of
$\psi (q)$- extended quantum oscillator algebras. With the choice $\psi
_n(q)=\left[ R(q^n)!\right] ^{-1}$and $R(x)=\frac{1-x}{1-q}$ we arrive at
the well known $q$-deformed oscillator. Because the reduced incidence
algebra $R(L(S))$ is isomorphic to the algebra $\Phi _\psi $ of $\psi $
-exponential formal power series - the $\psi $-extensions of finite
operator calculus provide a vast family of representations of $R(L(S))$.
\end{abstract}
\small{
}

\section{Introduction}

The main aim of this paper is a presentation of $\psi $- extension of Rota`s
finite operator calculus simultaneously with an indication that this is a
natural and elementary method for formulation and treatment of $q$-extended
and posibly $R$-extended or $\psi $-extended models for quantum-like $\psi $%
-deformed oscillators in order to describe \cite{3} eventual processes with
parastatistical behavior.

We owe such $\psi $-extension in operator form to Viskov \cite{1,2}.
This is achieved by considering not only
polynomial sequences of binomial type but also of $\{s_n\}_{n\geq 1}$%
-binomial type where $\{s_n\}_{n\geq 1}$-binomiality is defined with help of
the generalized factorial $n_s!=s_1s_2s_3...s_n$ where $S=\{s_n\}_{n\geq 1}$
is an arbitrary sequence with the condition $s_n\neq 0$, $n\in N$.

A year after Viskov's paper \cite{2} and few years before Roman \cite{4}-%
\cite{8} Cigler and Kirchenhofer in \cite{9a,9b} set up
foundations of $q$-umbral calculus as an extension of the umbral calculus
defined in \cite{10} in terms of the algebra of linear functionals $P^{*}$ -
see also \cite{11} two years after \cite{9a,9b}.

In our presentation here we use the operator formulation as in Rota-Mullin
calculus \cite{12}. The $\psi $- extension Rota`s finite operator calculus
extends the content of Rota`s devise ''much is the iteration of the few'' -
much of the properties of special polynomials is the application of few basic
principles.

Apart from pure mathematical interest there exist also other motivations. We
mention here only just one based on the simple observation. Namely
$q$-oscillator algebras generators are the so called $\partial _q$-delta
operators $ Q(\partial _q)$ and their duals and these are the basic objects
of the $q$-extended finite operator calculus of Rota to be formulated.
(Of course $\partial _q \hat x - q\hat x\partial_{q}=id$.)

These $q$-oscillator algebras generators are encountered explicitly or
implicitly in \cite{13,14} and in many other
subsequent references - see \cite{15,16,17a,17b,17c} and references therein.
There $q$-Laguerre and $q$-Hermite polynomials appear \cite{17a,17b,17c} which
are just $\partial _\psi $-basic polynomial
sequences of the $\partial _\psi $-delta operators $Q(\partial _{\psi})$ for
$\psi _n(q)=\frac 1{R(q^n)!};$ $R(x)=\frac{1-x}{1-q}$ and corresponding
choice of $Q(\partial _{\psi})$ functions of $\partial _\psi $- see next
sections. The case $\psi _n(q)=\frac 1{R(q^n)!}$ : $n_\psi =n_R;$ $\partial
_\psi =\partial _R$ and $n_{\psi (q)}=n_{R(q)}=R(q^n)$ appears implicitly
in \cite{18} where advanced theory of general quantum coherent states is
beeing developed. In \cite{19} it was noticed that commutation relations for
the $q$-oscillator algebras generators from \cite{13,14} and others
(see also \cite{17a,17b,17c,16} ) in
appropriate operator variables might be given the form \cite{19}:%
$$
AA^{+}-\mu A^{+}A=1;\textrm{ }\mu =q^2
$$
if appropriate operator variables are chosen \cite{19}. In this connection
let us note that various $q$-deformations of the natural number $n$
for the Fock space representation of normalized eigenstates of $\mid n>$%
of excitation number operator $N$ are widely used in literature on quantum
groups and at least some families of quantum groups may be constructed from
$q$-analogues of Heisenberg-Weyl algebra (see for example:
\cite{13,14,19,20,21a,21b} ).
Note also that $q$-analogues of Heisenberg-Weyl algebra appear as special
cases of general quantum Clifford algebras - understood as Chevalley-Kahler
deformations of braided exterior algebras \cite{22a,22b,22c}.

The known important fact is that the $q$-commutation relation $AA^{+}-A^{+}A
=1$ leads (see \cite{16}) to the $q$-deformed spectrum of excitation number
operator $N$ and to various parastatistics \cite{3} . More possibilities
result from considerations of Wigner \cite{23} extended by the authors of
\cite{3}. We therefore hope that the $\psi (q)$-calculus of Rota to be
developed here might be useful in a $C^{*}$ algebraic \cite{16} description
of ''$\psi (q)$-quantum processes'' - if any - with various parastatistics
\cite{3}.

Note also that the ''usual'' commutation relation
$AA^{+}-qA^{+}A=1$,
known since middle of nineteenth century \cite{24} is very useful in such
problems as derivation of Bernoulli -Taylor formula with the rest of Cauchy
form, inversion of formal power series and Lagrange formula \cite{25a,25b}.
The $\psi$-counterparts are expected.

The second purpose of the introduction is to establish notations to be
subsequently used . Here and once for all $P\equiv${\bf F}[x] denotes the
algebra of polynomials over the field $F$, char $F=0$. The standard
$q$-deformation of the real number is of the form $x_q\equiv \frac{1-q^x}{1-q}$
and the Jackson`s $\partial _q$- derivative (see: \cite{26a}-\cite{27c} and
\cite{28}) is defined by
$$
(\partial _q\varphi )(x)=\frac{\varphi (x)-\varphi (qx)}{(1-q)x}
$$
We use the following notation for $q$-factorial
$$
n_q!=n_q(n-1)_q!;\;1_q!=0_q!=1;\;n_q!\xrightarrow[q-1]{}n!.
$$
As for the further notation we enclose now {\bf the main abbreviations} used
in this paper. At first let $\Im $ be the set of function sequences such that
$$\Im = \{\psi ;R\supset[a,b];q\in[a,b];\psi (q):Z\rightarrow F;\psi _0(q)=1;
\psi _n(q)\neq 0;\psi_{-n}(q)=0;n\in N\}.$$

We introduce the notations to be used in the sequel with help of $\partial
_{\psi}$ and $n_{\psi}$ symbols where $n_\psi \equiv \psi _{n-1}(q)\psi
_n^{-1}(q);$ $n_\psi !\equiv \psi _n^{-1}(q)\equiv n_\psi (n-1)_\psi
(n-2)_\psi (n-3)_\psi ...2_\psi 1_\psi $ and $0_\psi !=1$.
The symbol $\partial _\psi$ is defined as follows. Let $\partial _\psi x^n=$
$n_\psi x^{n-1};n\geq 0$. The operator $\partial_{\psi}$ is then linearly
extended and we call it the $\psi $-derivative. There are some special cases
on the way.
For $\psi _n(q)=\frac 1{R(q^n)!}:n_\psi =n_R;$ $\partial _\psi =\partial _R$
and $n_{\psi (q)}=n_{R(q)}=R(q^n)$ \cite{18}. If $R(x)=\frac{1-x}{1-q}$ then
$n_\psi=n_q$ and $\partial _R=\partial _q$ \cite{9a,9b,10}.
Using the above abbreviations we are able now to introduce the next ones.\\

Abbreviations
\begin{enumerate}
\item $(x+_{\psi} a)^{n}=\sum \limits_{k\geq 0}\binom{n}{k}_{\psi}%
a^kx^{n-k}\;$ where\;
$\binom{n}{k}_{\psi} \equiv \frac{n_\psi ^{\underline{k}}}{k_\psi !}$
\;and\;\\$n_\psi^{\underline{k}}=
n_\psi (n-1)_\psi (n-2)_\psi ...(n-k+1)_\psi $; $\binom{n}{k}_q$
is the Gauss polynomial in the variable $q$ with integer coefficients and
its symbol $\binom{n}{k}_q$ shares \cite{28} main properties typical for
$\binom{n}{k}$.

\item $(x+_\psi a)^n\equiv E^a(\partial _\psi )x^n;$ $E^a(\partial _\psi
)= \sum\limits_{n\geq 0} \frac{a^n}{n_\psi !}\partial _\psi ^n;$
$E^a(\partial _\psi)f(x)=f(x+_\psi a);$

Note however that $(x+_\psi a)^n\neq (x+_\psi a)^{n-1}(x+_\psi a).$

For $q=1$ $E^a\left( \frac d{dx}\right) $is called the shift operator in
\cite{12,29} while in \cite{10} it is to be the translation
operator. In \cite{9a,9b} $ E^a(\partial _\psi )$ is also named
translation operator. The operators
$E^y(Q)= \sum\limits_{n\geq 0}p_n(y)\frac{Q^n}{n_\psi !}$ are generalised
translation operators \cite{30}.

\item $\sum_\psi $ denotes the algebra of $\partial _\psi $- shift invariant
operators where a linear operator $T_{\partial _\psi }:P\rightarrow P$ is $
\partial _\psi $-shift invariant iff $\forall _{\alpha \in F}\left[
T_{\partial _\psi },E^a(\partial _\psi )\right] =0.$

\item $Q(\partial _\psi )$ denotes the $\partial _\psi $-delta operator.
( It is a specific formal series in $\partial _\psi $).

\item The symbol $x_{Q(\partial_{\psi })}$ denotes the operator dual to the
$Q(\partial _\psi )$ operator.\\
The map $x_{Q(\partial_{\psi })}$ is called an ''umbral shift operator'' in
\cite{10} and in the functional formulation of undeformed umbral calculus
\cite{10} it is the operator
adjoint to a derivation of $P^*$ the linear space of linear functionals on
$P$ - see: Theorem 5 in \cite{10} and see also 1.1.16 in \cite{9b}.
Of course : for $Q = id$ we have : $x_{Q(\partial_{R})}\equiv
x_{\partial_{R}}\equiv \hat x$. (see Definition \ref{defnfiveone})

\item The linear operator $\hat x_{\psi}$ is defined accordingly:
$\hat x_{\psi}:P\rightarrow P$; $\hat x_{\psi} x^n=
\frac{\psi _{n+1}(q)(n+1)}{\psi _{n}(q)} x^{n+1} =
\frac{n+1}{(n+1)_{\psi}}x^{n+1}$; $n\geq 0$. In special cases we write:
$\hat x_{R}:P\rightarrow P$,
$\hat x_{R}x^n=\frac{(n+1)}{R(q^{n+1})}x^{n+1}$, $n\geq 0$;
$\hat x_{q}:P\rightarrow P$, $\hat x_{q}x^n=\frac{n+1}{(n+1)_{q}}x^{n+1}$,
$n\geq 0$.

\item The Pincherle $\psi $-derivative is the linear map {\bf '}:
$\sum_\psi \rightarrow \sum_\psi $ defined by the commutator
$$T_{\partial_{\psi } }\textrm{{\bf '}}= T_{\partial_{\psi } }\hat x_{\psi }
-\hat x_{\psi}T_{\partial_{\psi } }
\equiv [T_{\partial_{\psi } },\hat x_{\psi }].$$

\item Let now $\{p_{n}\}_{n\geq 0}$ be the $\partial_{\psi}$-basic polynomial
sequence of the $\partial_{\psi}$-delta operator $Q(\partial_{\psi})=Q$;
then the $\hat q_{\psi ,Q}$- operator is a liner map
$\hat q_{\psi ,Q}:P\rightarrow P$ defined in the basis $\{p_{n}\}_{n\geq 0}$
by $\hat q_{\psi ,Q}p_{n}=
\frac{(n+1)_{\psi}-1}{n_{\psi}}p_{n}$; $n\geq 0$.\\
We call the $\hat q_{\psi ,Q}$ operator the
$\hat q_{\psi ,Q}$-{\em mutator operator}.
For $\psi _{n}(q)=\frac{1}{R(q^n)!}$, $\hat q_{R,id}x^n=\frac{R(q^{n+1})-1}
{R(q^n)}x^n$, $n\geq 0$. If in addition $R(x)=\frac{1-x}{1-q}$ then
$\hat q_{R,id}x^n=qx^n$.

\item The $\hat q_{\psi ,Q}$-mutator of $A$ and $B$ operators reads:
$AB - \hat q_{\psi ,Q}BA\equiv [A,B]_{\hat q_{\psi ,Q}}$.
\end{enumerate}

The organization of the paper is the following.\\

In section 2 we observe that the reduced incidence algebra $R(L(S))$ is
isomorphic to the algebra $\Phi_{\psi}$ of formal $\psi$-exponential power
series. This makes the link of what follows after - to combinatorics.

Then in section 3 we note - as expected - that the algebra $\Phi_{\psi}$
of $\psi$-exponential formal power series is isomorphic to the algebra
$\sum _{\psi}$ of $\partial_{\psi}$-shift invariant operators which is the
basement of $\psi (q)$-extended finite operator calculus. We use this
opportunity to introduce there primary objects of the calculus.

In section 4 we develop further this $\psi (q)$- calculus
providing few examples of application of principal theorems.

With section 5 we finish our considerations by observing \cite{9a,9b}
that one
can formulate $q$-extended finite operator calculus with help of ''quantum
$q$-plane'' $q$-commuting variables $A,B: AB - qBA\equiv [A,B]_{q} = 0$.

The question whether one may formulate the $\psi$-extended finite operator
calculus with help of a ''quantum $\psi$-plane''
$\hat q_{\psi ,Q}$-commuting variables $A,B: AB - \hat q_{\psi ,Q}BA\equiv
[A,B]_{\hat q_{\psi ,Q}} = 0$ is discussed there also.


\section{Incidence algebras - primary information and possibility of $\psi$%
-extensions}

Pierre Simon Laplace (1749-1827) introduced the correspondence between
{\em operations} on sets and on {\em operations} on formal power series.
This has been afterwards developed into the vast activity domain of nowadays`
"generatinfunctionology" \cite{31}. One of few turning points in this area was
the use of incidence algebras by Rota and his collaborators
\cite{12}.
These incidence algebras were invented as a still more
systematic technique for setting up wider class of generating functions
algebras encompassing more than classical algebras of ordinary, exponential,
Dirichlet, Eulerian etc. generating functions. Incidence algebras were
independently introduced also by Scheid \cite{32} and Smith
\cite{33a,33b,33c} - see also \cite{34}.

In this section we follow the Rota`s way of presentation of
finite operator calculus as its extensions admit the similar treatment.
In order to establish the notation in a selfcontained way and to express our
main starting observation let us recall the required notions.

\begin{defn} \label{2.1}
Let $I(P,{\bf F})=\{f; f: P\times P\rightarrow {\bf F}; f(x,y)=0;
\;unless \;x\pi y;\\ x,y\in P\}$ where {\bf $F$} is a field;
$char\,F=0$ and $(P,{\bf \pi})$ is locally finite partially ordered
set. Then $(I(P,{\bf F}), {\bf F}; +;\textrm{ }^{\bf *};\textrm{ }^{\bf
\circ} )$ is called the incidence algebra of $P$, where "${\bf +}$" and
"$\textrm{}^{\bf \circ} $" denote the sum of functions and usual
multiplication by scalars, while for any $f,g\in I(P,{\bf F})$ the following
product is defined\\
$$(f^*g)(x,y)=\sum_{z\in P}f(x,z)g(z,y).$$
\end{defn}

Let as recall that a partially ordered set is locally finite iff every its
segment $[x,y]$ is finite, where $[x,y]=\{z\in P; x\leq z\leq y\}$ hence the
summation above is finite as it is over segment $[x,y]$. Here now are some
examples \cite{12} for illustration purpose.

\begin{ex}{\em
Let $P$ be the set of nonnegative integers \\ $P=\{0,1,2,3,4,5,6,7,...\}$
and let ${\bf \pi} \equiv \leq$ then
$$I(P,{\bf F})=\{(a_{ij});a_{ij}=0\;i<j\}\subset
M_{\infty}({\bf F})$$
i.e. $(I(P,{\bf F}), {\bf F}; {\bf +};\textrm{ }^{\bf *};\textrm{ }^
{\bf \circ})$ is represented by the algebra of upper triangular infinite
matrices over field ${\bf F}$.}
\end{ex}

\begin{ex}{\em
The algebra of formal power series is isomorphic to incidence algebra $R(P)$
where $(P;{\bf \pi})\equiv (P;\leq)$ and $P\equiv {\bf N}\cup \{0\}$.
(As a matter of fact $R(P)$ is the so called standard \underline {reduced}
incidence algebra see - below \cite{12}).\\
The isomorphism mentioned above is given by the bijective correspondence
$\phi$
\[
\sum\limits_{n \ge 0} {a_{n} z^{n}} \buildrel {\varphi}  \over
\longrightarrow f \equiv \left\{ {f_{ij} ;f_{ij} = \left\{
{{\begin{array}{*{20}c}
 {a_{i - j} \quad i \le j:\quad i,j \in P} \\
 {0\quad \quad otherwise} \\
\end{array}} } \right.} \right\}
\]
where $h\equiv f^*g$ with $f,g,h\in R(P)$  corresponds to convolution of
$\phi^{-1}(f)$ and $\phi^{-1}(g)$  i.e.
$$
h(i,j)=\sum_{i\leq k\leq j}f(i,k)g(k,j)=
\sum_{i\leq k\leq j}a_{k-i}b_{j-k} \equiv \sum_{r=0}^{n}a_{r}b_{n-r}
$$
after setting $r=k-i$ and $n=j-i$.}
\end{ex}

\begin{ex}{\em
The algebra of formal exponential power series is isomorphic to incidence
algebra $R(L(S))$; where $L(S)=\{A; A\subseteq S; |A|<\infty \}$ and $S$ is
countable while $(L(S);\subseteq)$ is partially ordered set - ordered by
inclusion. As a matter of fact $R(L(S))$ is the so called reduced incidence
algebra of the partially ordered set (poset) $L(S)$.\\
The isomorphism mentioned above is given by the bijective correspondence
$\phi$:
\[
F\left( {z} \right) \equiv \sum\limits_{n \ge 0} {\frac{{a_{n}
}}{{n!}}z^{n}\buildrel {\varphi}  \over \longrightarrow f = \left\{ {f\left(
{A,B} \right) = \left\{ {{\begin{array}{*{20}c}
 {a_{\left| {B - A} \right|} \quad A \subseteq B} \hfill \\
 {0\quad \quad otherwise} \hfill \\
\end{array}} ;\quad A,B \in L\left( {S} \right)} \right.} \right\}}
\]
where the product $h=f^*g$ with $f,g,h\in R(L(S))$  corresponds to binomial
convolution of  $F\equiv \phi ^{-1}(f)$ and $\phi ^{-1}(g)\equiv G$
i.e.  for $H(z)\equiv \sum\limits_{n\geq 0}
\frac{c_n}{n!}z^n$ ($H\equiv \phi ^{-1}(h)$) and
$G(z)\equiv \sum\limits_{n\geq 0}\frac{b_n}{n!}z^n$ we get $c_n=
\sum\limits_{k\geq 0}^{n}\binom{n}{k}a_{k}b_{n-k}$.\\
With "at the point" convergence $I(P;{\bf F})$ becomes a topological algebra
.}
\end{ex}

The $q$-extension or $\psi$-extension of the above examples and definitions
is automatic. It amounts to changes in notation
$n\rightarrow n_q\rightarrow n_{\psi}$
where - recall it again - $n_{\psi}\equiv \psi _{n-1}(q)\psi_{n}^{-1}(q)$;\\
$n_{\psi}! \equiv \psi_{n}^{-1}(q) \equiv
n_{\psi}(n-1)_{\psi}(n-2)_{\psi}(n-3)_{\psi}...2_{\psi}1_{\psi}$ and
$0_{\psi}!=1$ while $\psi \in \Im$ where\\
$\Im =\{\psi ;R\supset
[a,b];q \in [a,b];\psi (q):Z\rightarrow F;\psi _0(q)=1;\psi _n(q)\neq 0;\psi
_{-n}(q)=0;n\in N\}.$\\
Incidence algebras do characterise their correspondent partially ordered sets
as follows.

\begin{thm}
Let $P$, $Q$ be locally finite partially ordered sets. Let $I(P;{\bf F})$
and $I(Q;{\bf F})$ algebras be isomorphic. Then $P$ and $Q$ are isomorphic.
\end{thm}

Reduced incidence algebras and incidence coefficients are of the more
frequent use where the reduced incidence algebras are obtained as quotients
of incident algebras segments` families and an order compatible equivalence
relation. They corresponds to formal series of various kind. The incident
coefficients in their turn are generalisation of the binomial coefficients
\cite{12}.

\begin{defn}
Let ${\bf \sim}$ denote an equivalence relation defined on the family $S(P)$
of segments of $P$ - a locally finite partially ordered set.
Let $f,g\in I(P;{\bf F})$ be such that for $[x,y], [u,v]\in S(P)$ and
$[x,y] \sim [u,v]$  equalities $f(x,y)=f(u,v)$ and $g(x,y)=g(u,v)$ take place.
If $(f^*g)(x,y)=(f^*g)(u,v)$ $\forall$ $[x,y],[u,v]$ such that $[x,y] \sim
[u,v]$ then the relation "$\sim$" is said to be order compatible.
\end{defn}

For detailed information on properties and basic facts about incidence
algebras and compatible equivalence relation on locally finite partially
ordered sets see \cite{34,12,33a,33b,33c}.

\begin{defn}
Let $P$ be a locally finite partially ordered set equipped with a
compatible equivalence relation $\sim$ on $S(P)$. The set of all functions
defined on $S(P)/ \sim$ with the product defined below
is called the reduced incidence algebra $R(P;\sim)$.
\end{defn}

In order to define the product of  $f: S(P)/\sim \rightarrow {\bf F}$ and
$g: S(P)/\sim \rightarrow {\bf F}$ referred to in the definition above let us
denote by $\alpha ,\beta ,...$  the nonempty equivalence classes of segments
of $P$ i.e. $\alpha ,\beta ,...\in S(P)/\sim$ and let us call them \cite{12}
{\em types}.

\begin{defn}
$(Map(S(P)/{\bf \sim}; {\bf F}),{\bf F};
{\bf +};\textrm{ }^{\bf *};^{\bf \circ}) \equiv {\bf R(P;{\bf \sim})}$ is an
algebra under the multiplication "$\textrm{ }^*$" defined as follows:
$$Map(S(P)/{\bf \sim}; {\bf F})\ni f,g\rightarrow h:=f^*g;$$
$$S(P)/\sim \ni \alpha \rightarrow h(\alpha):=\sum_{(\Lambda)}
\genfrac{[}{]}{0pt}{}{\alpha}{\beta \gamma}f(\beta)g(\gamma),$$
where the sum $\sum\limits_{(\Lambda)}$ ranges over all ordered pairs
$(\beta ,\gamma)$ of all types and the brackets
$\genfrac{[}{]}{0pt}{}{\alpha}{\beta \gamma}$ are defined below.
\end{defn}

\begin{defn}
$\genfrac{[}{]}{0pt}{}{\alpha}{\beta \gamma}:=$ the number of such distinct
elements $z$ in a segment $[x,y]$ of type $\alpha$ and such elements $z$
that $[x,z]$ is of type $\beta$ while $[z,y]$ is of type $\gamma$.
\end{defn}

One may prove \cite{12} that the reduced incidence algebra $R(P;{\bf \sim})$
$\{$i.e. the incidence algebra modulo ${\bf \sim}\}$ is isomorphic to a
subalgebra of the incidence algebra of $P$.\\
Now we may formulate the main observation of this section which links reduced
incidence algebras with their representations by $\psi$-{\em extended}
Rotas` calculus (see: next sections).\\ For that to see it now and then to
explore let us announce ahead the existence (theorem~\ref{ththreethree}) of
the algebra isomorphism $\sum_{\psi} \approx \Phi_{\psi}$ where $\sum_{\psi}$
denotes the algebra of $\partial_{\psi}$ -shift invariant operators and
$\Phi_{\psi}$ is the algebra of  $\psi$-exponential formal power series. This
is to be combined with the following observation.\\
{\bf Starting Observation}\\

The algebra $\Phi_\psi$ of $\psi$-exponential formal power series is
isomorphic to the reduced incidence algebra $R(L (S))$. The isomorphism
$\phi$ is given by the bijective correspondence
$$
F_{\psi}(z)=\sum_{n\geq 0}\frac{a_n}{n_{\psi}!}z^n\xrightarrow{\phi}
f=\left\{{}f(A,B)= \begin{cases} a_{|A-B|}&\,A\leq B\\0&\textrm{otherwise}
\end{cases};A,B\in L(S)\right\}
$$
Here $h:=f^*g$ - where $f,g,h\in (R(L(S))$ - corresponds to $\psi$-binomial
convolution i.e. for
$$
H_{\psi}(z)=\sum_{n\geq 0}\frac{c_n}{n_{\psi}!},\;
G_{\psi}(z)=\sum_{n\geq 0}\frac{b_n}{n_{\psi}!}z^{n}\; and \;\,
[n]h(z)\equiv [n](f^*g)(z)\equiv c_n;
$$
$$c_n=\sum_{k\geq 0}^n\binom{n}{k}_{\psi}a_{k}b_{n-k}.$$

\begin{rem}{\em
Before coming over to presentation of $\psi$-extension of finite operator
calculus let us remark that variety of enumerative problems are instances
\cite{34} of the general problem of inverting indefinite sums \cite{35}
ranging over a partially ordered set \cite{34}.

The inversion can be carried out with the use of an analogue of difference
operator (relative to this partial order).

Next - the indefinite sums ranging in the locally finite partially ordered
set $P$ are analogues of indefinite integral - while various difference
operators (for example $\partial_{\psi}$-delta operators) are analogues of
differentiation operator $D$.\\
Namely let $f:P\rightarrow F$ be a function defined on a locally finite
partially ordered set $P$. Let the indefinite sum $g$ be calculated via
$g(x)=\sum\limits_{y\leq x}f(y)$ then the differentiation is being
realized by $f(x)=\sum\limits_{y\leq x}g(y)\mu (y,x)$ i.e. by the M\"obius
inversion formula \cite{34}, where $\mu \in I(P,F)$ is the M\"obius function
on the locally finite partially ordered set $P$ inverse to zeta function
$\zeta \in I(P,F)$ determined by $\zeta(x,y)=1\; for\; x\leq y$ and
$\zeta(x,y)=0\;otherwise$\; where of course $x,y\in P$.}
\end{rem}

In the case of $\Phi_{\psi}$ representation of the reduced incidence algebra
$R(L(S))$ the incidence coefficients correspond to $\psi$-binomial
coefficients and the types correspond to $\psi$-extended integers.

It is known that the binomiality of basic polynomial sequences resulting
through isomorphisms
$\sum _{\psi} \approx \Phi_{\psi} \approx R(L(S))$ has
transparent combinatorial interpretation (see: \cite{36,37} and
\cite{31} - page 91. It remains to be an open question for the present author
what is the role and scope of possible applications to combinatorics of
$\psi$-binomiality.

Anyhow as we shall see it in the sequel - at least special $\psi$-extensions
supply a mathematical underpinning for $\psi(q)$-deformed quantum oscillator
algebras .

\section{Finite Operator $\psi$-Calculus - an Introduction}

In sections 3 and 4 we give an exposition of the beginnings of $\psi
$-extension of umbral calculus. This includes among others proofs of The
First and The Second Expansion Theorems as well as Rodrigues formula and
Spectral Theorem .

Let us however start from the beginning i.e. from the idea.\\
\textit{The objective} of \cite{12} was a unified theory of special polynomials.
We extend this objective to encompass also correspondent $\psi $-extended
families of polynomials i.e. as a matter of fact - to encompass all those
polynomial sequences $\left\{ {s_{n} \left( {x} \right)} \right\}_{n =
0}^{\infty}  $ , deg s$_{n} $= n. (E. Loeb \cite{38}) which may be considered as a
Sheffer polynomial sequence with respect the corresponding $\partial _{\psi
} $-delta operator $Q\left( {\partial _{\psi} }  \right)$.
A still further generalisation is to be found in Markowsky paper \cite{29}
where generalised Sheffer polynomial sequences are just polynomial sequences
$\left\{ {s_{n} \left( {x} \right)} \right\}_{n = 0}^{\infty}  $; deg s$_{n
}= n$  and the role of $\partial _{\psi}  $-delta operators
$Q\left( {\partial _{\psi} }  \right)$ is taken over by generalised
differential operators \cite{29}. Many of the results of $\psi $-calculus to
be developed here may be extended to Markowsky $Q$-umbral calculus where
$Q$ stands for a generalised difference operator i.e. the one
lowering the degree of any polynomial by one.\\
\textit{The way} to achieve our $\psi$-goal as in \cite{12} is to exploit
the duality between the $\hat {x}$ and $\frac{{d}}{{dx}}$ the
predecessors of the delta operator notion and its dual.\\
\textit{The technique} used and co-invented mostly by \cite{12} is of the past
century origin and it is the so called symbolic calculus.

In this section we shall refer all the time to Rota \cite{12} and we shall
follow his way of presentation. The language as well as notation have been
designed to reflect this resolution. Let us then start with recalling the
basic definitions.

\begin{defn}
Let $E^{y}\left( {\partial _{\psi} }  \right) \equiv exp_{\psi}  \{
y\partial _{\psi}  \} = \sum\limits_{k = 0}^{\infty}  {\frac{{y^{k}\partial
_{\psi}  ^{k}}}{{n_{\psi}  !}}} $ be the linear operator $E^{y}\left(
{\partial _{\psi} }  \right)$:\textit{ P$ \to $ P}.
Let $T_{\partial _{\psi} }  \;\;:\;P \to P$ be a linear operator; then
$T_{\partial _{\psi} }  \;$is $\partial _{\psi} $-shift invariant iff
$$\forall \alpha \in F; \; \left[ T_{\partial _{\psi}},E^{\alpha} \left(
{\partial _{\psi}}
\right) \right] = 0.$$
\end{defn}
\textbf{Notation}: $\sum_{\psi} $denotes the algebra of
\textit{F}-linear $\partial _{\psi}  $-shift invariant operators
$T_{\partial _{\psi} }$.

\begin{defn}
Let $Q\left( {\partial _{\psi} }  \right):P \to P$ be a linear operator
$Q\left( {\partial _{\psi} }  \right)$ is a $\partial _{\psi}$-delta operator
iff

\begin{enumerate}
\renewcommand{\labelenumi}{(\alph{enumi})}
\item $Q\left( {\partial _{\psi} }  \right)$ is $\partial _{\psi}  $-shift
invariant;
\item $Q\left( {\partial _{\psi} }  \right)\left( {id} \right) =
const \ne 0.$
\end{enumerate}
\end{defn}

\begin{obs}\label{threeone}
Let $Q\left( {\partial _{\psi} }\right) \in
\sum_{\psi}$ be the $\partial _{\psi}$-delta operator. Then for
every constant polynomial $a \in P$ we have
 $Q\left( {\partial _{\psi} }  \right) a  = 0$.

\begin{proof}
Recall that $\forall \;\alpha \in F;\quad \left[ {Q\left( {\partial _{\psi}
} \right),E^{\alpha} \left( {\partial _{\psi} }  \right)} \right] = 0$ then
by linearity of $Q\left( {\partial _{\psi} }  \right)$ we have\\
($Q\left( {\partial _{\psi} }  \right)E^{a}\left( {\partial _{\psi} }
\right)$)(x) = $Q\left( {\partial _{\psi} }  \right)$ (x+a) = $Q\left(
{\partial _{\psi} }  \right)$ (x) + $Q\left( {\partial _{\psi} }  \right)$
(a) = c + $Q\left( {\partial _{\psi} }  \right)$ (a)\\
and at the same time ($E^{a}\left( {\partial _{\psi} }  \right)Q\left(
{\partial _{\psi} }  \right)$)(x) = $E^{a}\left( {\partial _{\psi} }
\right)$ (c) = c.
\end{proof}
\end{obs}

\begin{obs}
If $p \in P$, deg p = n then $deg\left( {Q\left( {\partial
_{\psi} }  \right)p_{n}}  \right) = n - 1$.

\begin{proof} The proof goes like in \cite{12}. We ``just replace'' shift
invariance by $\partial _{\psi}  $-shift invariance .\\
 $Q\left( {\partial _{\psi} }  \right)$ ($E^{a}\left( {\partial _{\psi} }
\right)$ (x$^{n}$)) = $Q\left( {\partial _{\psi} }  \right)$ (x +$\psi $
a)$^{n} \quad  \equiv \sum\limits_{k \ge 0} {\left( {{\begin{array}{*{20}c}
 {n} \hfill \\
 {k} \hfill \\
\end{array}} } \right)} _{\psi}  a^{k}Q\left( {\partial _{\psi} }
\right)x^{n - k}$=$E^{a}\left( {\partial _{\psi} }  \right)$ ( $Q\left(
{\partial _{\psi} }  \right)$ (x$^{n}$)) $ \equiv r\left( {x + _{\psi}  a}
\right)$ , where $\left( {{\begin{array}{*{20}c}
 {n} \hfill \\
 {k} \hfill \\
\end{array}} } \right)_{\psi}  \equiv \frac{{n_{\psi} ^{\underline {k}}
}}{{k_{\psi}  !}}\,$.\\
Hence $r\left( {a} \right) = \sum\limits_{k \ge 0}
{\left( {{\begin{array}{*{20}c}
 {n} \hfill \\
 {k} \hfill \\
\end{array}} } \right)} _{\psi}  a^{k}Q\left( {\partial _{\psi} }
\right)x^{n - k}\left| {_{x = 0}}  \right.$ . The coefficient of $a^{n}$ is
therefore equal to zero according to the observation~\ref{threeone}. At the
same time the coefficient of $a^{n - 1}$ is equal to
$\left( {{\begin{array}{*{20}c}
 {n} \hfill \\
 {n - 1} \hfill \\
\end{array}} } \right)_{\psi}  \equiv c \ne 0$
\end{proof}
\end{obs}

\begin{defn}
A polynomial sequence $\left\{ {p_{n}}  \right\}_{n \ge 0} $, deg
$p_{n} = n$ such that
\begin{enumerate}
\renewcommand{\labelenumi}{(\alph{enumi})}
\item $p_{o} \left( {x} \right) = 1$;
\item $p_{n} \left( {0} \right) = 0$; $n > 0$;
\item $Q\left( {\partial _{\psi} }  \right)p_{n} = n_{\psi}
p_{n - 1}$
\end{enumerate}
is called the $\partial _{\psi}  $-basic polynomial sequence of the
$\partial _{\psi}  $-delta operator $Q\left( {\partial _{\psi} }  \right)$.
\end{defn}

The condition b) $p_{n} \left( {0} \right) = 0$, $n > 0$ used in \cite{12}
is supeflouous as from a) and c) $\psi $-binomiality is
easily proven. From the $\psi $-binomiality property one derives b).

\begin{prop}
Every $\partial _{\psi}  $-delta operator
$Q\left( {\partial _{\psi} }  \right)$ has the unique sequence of $\partial
_{\psi}  $-basic polynomials i.e.
\[
Q\left( {\partial _{\psi} }  \right)\;{\begin{array}{*{20}c}
 {\buildrel {1:1} \over \longrightarrow}  \hfill \\
 {\mathrel{\mathop{\kern0pt\longleftarrow}\limits_{1:1}}}  \hfill \\
\end{array}} \;\left\{ {p_{k}}  \right\}_{0}^{\infty}
\]
\begin{proof}
For $n=0$ put $p_{o} \left( {z} \right) = 1$, for $n=1$ put $p_{1} \left(
{x} \right) = \frac{{x}}{{Q\left( {\partial _{\psi} }  \right)\left( {id}
\right)}}$. Then inducing on \textit{n} assume that
\{\textit{p}$_{k}$(\textit{z})\} have been defined for $k<n$. From this
inductive assumption we infer that \textit{p}$_{n} $ is defined uniquely.
For that to see it is enough to notice that for any \textit{p$ \in $ P}, deg
\textit{p} = \textit{n} i.e.
 $$p\left( {z} \right) = az^{n} + \sum\limits_{k = 0}^{n - 1} {c_{k} p_{k}
\left( {z} \right)} \textrm{and} a \ne 0$$
we have $Q\left( {\partial _{\psi} }  \right)p\left( {z} \right) = aQ\left(
{\partial _{\psi} }  \right)x^{n} + \sum\limits_{k = 1}^{n - 1} {c_{k} k_{q}
p_{k - 1} \left( {z} \right)} $  and  $degQ\left( {\partial _{\psi} }
\right)\left( {z^{n}} \right) = n - 1$.\\
Hence there exist a unique choice of constants c$_{1}$, ... ,c$_{n-1}$ for
which
$$Q\left( {\partial _{\psi} }  \right) p=n_{q} p_{n-1;}$$

This determines $p \equiv p_{n} $ uniquely except for the constant term
$c_{o}$ which is however determined uniquely by the condition $p_{n}(0)=0$,
$n>0$.
\end{proof}
\end{prop}

Inspired by the predecessors $\hat {x}$ and $\frac{{d}}{{dx}}$ of
the notions developed in \cite{12} we introduce the next basic notion.

\begin{defn}
A polynomial sequences $\left\{ {p_{n}}  \right\}_{o}^{\infty}  $ is of
$\psi $-binomial type if it satisfies the recurrence
$E^{y}\left( {\partial _{\psi} }  \right)
p_{n} \left( {x} \right)
 \equiv
p_{n} \left( {x + _{\psi}  y} \right) = \sum\limits_{k \ge 0} {\left(
{{\begin{array}{*{20}c}
 {n} \hfill \\
 {k} \hfill \\
\end{array}} } \right)} _{\psi}  p_{k} \left( {x} \right)p_{n - k} \left(
{y} \right)$ where $\left( {{\begin{array}{*{20}c}
 {n} \hfill \\
 {k} \hfill \\
\end{array}} } \right)_{\psi}  \equiv \frac{{n_{\psi} ^{\underline {k}}
}}{{k_{\psi}  !}}\,$.
\end{defn}

\begin{thm}\label{ththreeone}
 $\left\{ {p_{n}}  \right\}_{o}^{\infty}  $ is a $\partial _{\psi}  $-basic
sequence of some $\partial _{\psi}  $-delta operator $Q\left( {\partial
_{\psi} }  \right)$ iff it is a sequence of $\psi $-binomial type.
\end{thm}
\begin{proof}
$\textrm{   }$
\begin{enumerate}
\renewcommand{\labelenumi}{(\bf {\Alph{enumi}})}
\item $Q\left( {\partial _{\psi} }  \right)^{k}p_{n} \left( {x}
\right) = n_{\psi} ^{\underline {k}}  p_{n - k} \left( {x} \right)$
therefore for $k=n\textrm{ }$\\
$\left[ {Q\left( {\partial _{\psi} }  \right)^{n}p_{n}
\left( {x} \right)} \right]\left|{_{x = 0}} \right. = n_{\psi}!$ while for
$0<k<n\textrm{ }$ $\left[ {Q\left( {\partial _{\psi} }  \right)^{k}p_{n}
\left( {x} \right)} \right]\left|{_{x = 0}} \right. = 0$ hence
$p_{n} \left( {x} \right) = \sum\limits_{k \ge 0} {\frac{{p_{k} \left( {x}
\right)}}{{k_{\psi} !}}} \left[ {Q\left( {\partial _{\psi} }
\right)^{k}p_{n} \left( {x}\right)} \right]\left|{_{x = 0}} \right. $
and by linearity argument $\forall p \in P$
\[
p\left( {x} \right)=
\sum\limits_{k \ge 0} {\frac{{p_{k} \left( {x} \right)}}{{k_{\psi}  !}}}
\left[ {Q\left( {\partial _{\psi} }  \right)^{k}p\left( {x} \right)}
\right]\left|{_{x = 0}} \right.
\]
With the choice
$p\left( {x} \right) = p_{n} \left( {x + _{\psi}  y}
\right)$ where $y$ is a parameter - one gets\\
$p_{n} \left( {x + _{\psi}  y} \right)=\sum\limits_{k \ge 0}
{\frac{{p_{k} \left( {x} \right)}}{{k_{\psi}  !}}}
\left[ {Q\left( {\partial _{\psi} }  \right)^{k}p_{n} \left( {x + _{\psi}
y} \right)} \right]\left|{_{x = 0}} \right.$ from which one infers that\\
\begin{multline*}
\left[ {Q\left( {\partial _{\psi} }  \right)^{k}p_{n} \left( {x + _{\psi}
y} \right)} \right]\left|{_{x = 0}} \right. =
\left[ {Q\left( {\partial _{\psi} }  \right)^{k}E^{y}\left( {\partial _{\psi
}}  \right)p_{n} \left( {x} \right)} \right]\left|{_{x = 0}} \right. =\\
=\left[ {E^{y}\left( {\partial _{\psi} }  \right)Q\left( {\partial _{\psi} }
\right)^{k}p_{n} \left( {x} \right)} \right]\left|{_{x = 0}} \right.
=\left[ {E^{y}\left( {\partial _{\psi} }  \right)n_{\psi} ^{\underline {k}
} p_{n - k} \left( {x} \right)} \right]\left|{_{x = 0}} \right. =\\
= n_{\psi}^{\underline {k}}  p_{n - k} \left( {x + _{\psi}  y}
\right)\left|{_{x =0}} \right. = n_{\psi} ^{\underline {k}}  p_{n - k}
\left( {y} \right).
\end{multline*}
Altogether $$p\left( {x + _{\psi}  y} \right)=
\sum\limits_{k \ge 0} {\frac{{p_{k} \left( {x} \right)}}{{k_{\psi}  !}}}
n_{q}^{\underline {k}}
p_{n - k} \left( {y} \right)=
\sum\limits_{k \ge 0} {\left( {{\begin{array}{*{20}c}
 {n} \hfill \\
 {k} \hfill \\
\end{array}} } \right)} _{\psi}  p_{k} \left( {x} \right)p_{n - k} \left(
{y} \right) \eqno (*)$$

\item  suppose now that $\left\{ {p_{n}}  \right\}_{o}^{\infty}  $ is
a sequence of $\psi $-binomial type. Setting in (*) $y=0$ we get
\begin{multline*}
p_{n} \left( {x} \right) = \sum\limits_{k \ge 0} {\left(
{{\begin{array}{*{20}c}
 {n} \hfill \\
 {k} \hfill \\
\end{array}} } \right)} _{\psi}  p_{k} \left( {x} \right)p_{n - k} \left(
{0} \right) \quad =\\
\quad
p_{n} \left( {x} \right)
p_{n} \left( {0} \right)+
\quad
n_{\psi}
p_{n - 1} \left( {x} \right)
p_{1} \left( {0} \right) \quad +
\quad
\left( {{\begin{array}{*{20}c}
 {n} \hfill \\
 {2} \hfill \\
\end{array}} } \right)_{\psi}
p_{n - 2} \left( {x} \right)
p_{2} \left( {0} \right) \quad + \quad ...
\end{multline*}
from which we infer that $p_{0} \left( {x} \right)=1$ and $p_{n} \left( {0}
\right) = 0$ for $n>0$ . It is sufficient now to define the ${\partial
_{\psi}}$-delta operator $Q \left( {\partial _{\psi}} \right)$ corresponding
to $\left\{ {p_{n}}  \right\}_{o}^{\infty}$. We define it uniquely
according to:
\begin{enumerate}
\item $Q\left( {\partial _{\psi} }  \right)p_{0} \left( {x} \right) = 0$,

\item $Q\left( {\partial _{}}  \right)p_{n} = n_{\psi}  p_{n - 1} $ for $n>0$,

\item $Q\left( {\partial _{\psi} }  \right)$is linear.
\end{enumerate}
We now prove that $Q\left( {\partial _{\psi} }  \right)$ is $\partial _{\psi
} $-shift invariant. For that to do use
\[
p_{n} \left( {x + _{\psi}  y} \right)=
\sum\limits_{k \ge 0} {\frac{{p_{k} \left( {x} \right)}}{{k_{\psi}  !}}}
Q\left( {\partial _{\psi} }  \right)^{k}
p_{n} \left( {y} \right)
\]
which by linearity argument extends to any $p \in P$
\[
p\left( {x + _{\psi}  y} \right)=
\sum\limits_{k \ge 0} {\frac{{p_{k} \left( {x} \right)}}{{k_{\psi}  !}}}
Q\left( {\partial _{\psi} }  \right)^{k}
p\left( {y} \right).
\]
Now \cite{12} replace $p$ by $Q\left( {\partial _{\psi} }  \right)p$ and
interchange $x$ and $y$ thus getting
\[
\quad
\left( {Q\left( {\partial _{\psi} }  \right)p} \right)\left( {x + _{\psi}
y} \right) =
\sum\limits_{k \ge 0} {\frac{{p_{k} \left( {y} \right)}}{{k_{\psi}  !}}}
Q\left( {\partial _{\psi} }  \right)^{k + 1}
p\left( {x} \right).
\]
However note that\\
$\left( {Q\left( {\partial _{\psi} }  \right)p}
\right)\left( {x + _{\psi}  y} \right)$= $E^{y}\left( {\partial _{\psi} }
\right)\left( {Q\left( {\partial _{\psi} }  \right)p} \right)\left( {x}
\right)$=$E^{y}\left( {\partial _{\psi} }  \right)Q\left( {\partial _{\psi
}}  \right)p\left( {x} \right)$ and
\begin{multline*}
\sum\limits_{k \ge 0} {\frac{{p_{k} \left( {y} \right)}}{{k_{\psi}  !}}}
Q\left( {\partial _{\psi} }  \right)^{k + 1}
p\left( {x} \right) =
Q\left( {\partial _{\psi} }  \right)
\sum\limits_{k \ge 0} {\frac{{p_{k} \left( {y} \right)}}{{k_{q} !}}}
Q\left( {\partial _{\psi} }  \right)^{k}
p\left( {x} \right)=\\
=Q\left( {\partial _{\psi} }  \right)
\left( {p\left( {x + _{\psi}  y} \right)} \right)=
Q\left( {\partial _{\psi} }  \right)E^{y}\left( {\partial _{\psi} }
\right)p\left( {x} \right).\end{multline*}
\end{enumerate}
\end{proof}

\begin{thm} {\em (First $\psi$-Expansion Theorem)}\\
Let $T_{\partial _{\psi} } :P \to P$ be a $\partial _{\psi}  $-shift
invariant operator. Let $Q\left( {\partial _{\psi} }  \right)$ be a
$\partial _{\psi}  $-delta operator with $\partial _{\psi}  $-basic
polynomial sequence $\left\{ {p_{n}}  \right\}_{o}^{\infty}$. Then
$$T_{\partial _{\psi} }  = \sum\limits_{n \ge 0} {\frac{{a_{n}} }{{n_{\psi}
!}}} Q\left( {\partial _{\psi} ^{}}  \right)^{n} \quad ,\textrm{ where }
a_{k} = {\left[ {T_{\partial _{\psi} }  p_{k} \left( {z} \right)} \right]}
\left|{_{z = 0}} \right. .$$

\begin{proof} Due to Theorem \ref{ththreeone} the proof is the same as in
\cite{12}. Namely the $\psi$-binomial formula
\[
p_{n} \left( {x + _{\psi}  y} \right)=
\sum\limits_{k \ge 0} {\frac{{p_{k} \left( {x} \right)}}{{k_{\psi}  !}}}
Q\left( {\partial _{\psi} }  \right)
^{k}
p_{n} \left( {y} \right)
\]
acted upon by $\partial _{\psi}  $-shift invariant operator $T_{\partial
_{\psi} }$ extends by linearity argument to
\[
T_{\partial _{\psi} }  p\left( {x + _{\psi}  y} \right)=
\sum\limits_{k \ge 0} {\frac{{T_{\partial _{\psi} }  p_{k} \left( {x}
\right)}}{{k_{\psi}  !}}}
Q\left( {\partial _{\psi} }  \right)^{k}
p\left( {y} \right), \textrm{where} y \textrm{is treated as a parameter.}
\]
Setting above $x=0$ and exchanging mutually $x$ with $y$ symbols we arrive
at
\[
T_{\partial _{\psi} }  p\left( {x} \right)=
\sum\limits_{k \ge 0} \frac{\left[ T_{\partial _{\psi} }p_{k} \left(
{y} \right) \right]_{y = 0}}{{k_{\psi}!}}
Q\left( {\partial _{\psi} }  \right)^{k}
p\left( {x} \right).
\]
\end{proof}
\end{thm}

\begin{rem}{\em
The first expansion theorem might serve us to express popular
$q$-deformations of delta operators in terms of the others.\\ (Recall:
$n_{\psi}  \equiv \psi_{n - 1} \left( {q} \right)\psi _{n}^{ - 1}
\left( {q} \right)$. For $\psi_{n} \left( {q} \right) =
\frac{{1}}{{R\left( {q^{n}} \right)!}}$ :
$n_{\psi} = n_{R} $ ; $\partial _{\psi} = \partial _{R} $ and $n_{\psi
\left( {q} \right)} = n_{R\left( {q} \right)} = R\left( {q^{n}} \right)$
and if in addition $R\left( {x} \right) = \frac{{1 - x}}{{1 - q}}$ then
$n_{\psi} = n_{q} $ and $\partial _{R}  = \partial _{q} $ ). For
example the easy expansion of difference operator $\Delta = \sum\limits_{n
\ge 1} {\frac{{\delta _{n}} }{{n!}}} \frac{{d^{n}}}{{dx^{n}}}$ where $\delta
_{n} = \left[ {\Delta x^{n}} \right]_{x = 0} $ =1 suggests a natural
$q$-deformation of difference operator $\Delta $ in the form $\Delta _{q} =
\sum\limits_{n \ge 1} {\frac{{1_{q}} }{{n_{q} !}}} \partial _{q}^{n} $. Of
course every $\partial _{q} $-shift invariant operator $T_{\partial _{q}}
\;$is of the form $T_{\partial _{q}}  = \sum\limits_{n \ge 0} {\frac{{a_{n}
}}{{n_{q} !}}} \partial _{q}^{n} $ and in particular every $\partial _{q}
$\textit{-delta} operator $Q\left( {\partial _{q}}  \right)$ has an
expansion with help of any other one. For example for delta operators
$\frac{{d}}{{dx}}$ and $\Delta $ one has $\frac{{d}}{{dx}} = \sum\limits_{k
\ge 1} {\frac{{d_{k}} }{{k!}}\Delta ^{k}} $ where $d_{k} = \left[
{\frac{{d}}{{dx}}x^{\underline
 {k}} } \right]_{x = 0} = \left( { - 1}
\right)^{k - 1}\left( {k - 1} \right)!$ and correspondingly : $\partial _{q}
= \sum\limits_{k \ge 1} {\frac{{d_{k}^{\left( {q} \right)}} }{{k_{q}
!}}\Delta _{q} ^{k}} $ where $d_{k}^{\left( {q} \right)} $ might be
calculated using the explicit form of $\partial _{q} $-basic polynomial
sequence of the $\Delta _{q} $. This form is obtained with help of statement
(3) of the theorem~\ref{thfourone} - see examples at the end of section 4.}
\end{rem}

\begin{thm}\label{ththreethree}
$\Sigma_{\psi} \approx \Phi _{\psi}  $\\
Let $Q\left( {\partial _{\psi} }  \right)$ be a $\partial _{\psi}  $-delta
operator and let $\Phi _{\psi}  $ be \textit{the algebra of formal
exp}$_{\psi}  $\textit{ s}eries of t$ \in $\textit{F}\textbf{} over the same
field \textit{F} for which $Q\left( {\partial _{\psi} }  \right)$ is defined.
Then there exists an isomorphism $\varphi $:\\ $\varphi $ : $\Phi _{\psi}
 \to \Sigma_{\psi}$ of the algebra $\Phi _{\psi}$
onto the algebra $\Sigma _{\psi}  $ of $\partial _{\psi}  $-shift
invariant operators $T_{\partial _{\psi} }  \;$ which carries
\[
f_{\psi}  \left( {t} \right) = \sum\limits_{k \ge 0}
\frac{a_{k}t^{k}}{k_{\psi}!}
\buildrel {into} \over \longrightarrow T_{\partial_{\psi}} =
\sum\limits_{k \ge 0} \frac{a_{k}}{k_{\psi}!} \partial_{\psi}^{k}.
\]
\end{thm}
\begin{proof}
The proof goes like in \cite{12} without significant changes. Indeed, with
$\varphi $ being obviously linear and onto by the first expansion theorem,
it is enough to show that $\varphi $ preserves products.\\
In the algebra $\Phi _{\psi}  $ the product is given by the $\psi
$-binomial convolution i.e. we have
\[
\left( {\;\sum\limits_{k \ge 0} {\frac{{a_{k}} }{{k_{\psi}  !}}} x^{k}\;}
\right)
\left( {\;\sum\limits_{k \ge 0} {\frac{{b_{k}} }{{k_{\psi}  !}}} x^{k}\;}
\right)=
\left( {\;\sum\limits_{k \ge 0} {\frac{{c_{k}} }{{k_{\psi}  !}}} x^{k}\;}
\right)
\]
where
\[
c_{n} = \sum\limits_{k \ge 0} {\left( {{\begin{array}{*{20}c}
 {n} \hfill \\
 {k} \hfill \\
\end{array}} } \right)_{\psi}  a_{k} b_{n - k}}  .
\]
Therefore with
\[
\Phi _{\psi} \ni s\left( {x} \right) = \sum\limits_{k \ge 0} {\frac{{b_{k}
}}{{k_{\psi}  !}}} x^{k}
\quad
\buildrel {on} \over \longrightarrow
\quad
\sum\limits_{k \ge 0} {\frac{{b_{k}} }{{k_{\psi}  !}}} Q\left( {\partial
_{\psi} }  \right)^{k}\quad = \,S_{\partial _{\psi} }  \in \Sigma _{\psi}
\]
it is enough to show that
\[
\left[ {T_{\partial _{\psi} }  S_{\partial _{\psi} } ^{} p_{n} \left( {x}
\right)} \right]_{x = 0} \; = \;c_{n} = \sum\limits_{k \ge 0} {\left(
{{\begin{array}{*{20}c}
 {n} \hfill \\
 {k} \hfill \\
\end{array}} } \right)_{\psi}  a_{k} b_{n - k}}
\]
and this is the case because - due to $p_{n} \left( {0} \right) = 0$ for
$n>0$ \& $p_{0} \left( {x} \right) = 1$ we have
\[
\left[ {T_{\partial _{\psi} }  S_{\partial _{\psi} }  p_{n} \left( {x}
\right)} \right]_{x = 0} \;=
\quad
\left[ {\left( {\sum\limits_{k \ge 0} {\frac{{a_{k}} }{{k_{\psi}  !}}Q\left(
{\partial _{\psi} }  \right)^{k}\sum\limits_{r \ge 0} {\frac{{b_{r}
}}{{r_{\psi}  !}}Q\left( {\partial _{\psi} }  \right)^{r}}} }  \right)p_{n}
\left( {x} \right)} \right]_{x = 0} =
\]
\[
=\left[ {\sum\limits_{k \ge 0} {\sum\limits_{r \ge 0} {\frac{{a_{k}
}}{{k_{\psi}  !}}\frac{{b_{r}} }{{r_{\psi}  !}}Q\left( {\partial _{\psi} }
\right)^{k + r}p_{n} \left( {x} \right)}} }  \right]_{x = 0} =
\]
\[
=\left[ {\sum\limits_{k \ge 0} {\frac{{a_{k} b_{n - r}} }{{k_{\psi}
!\left( {n - r} \right)_{\psi}  !}}Q\left( {\partial _{\psi} }
\right)^{n}p_{n} \left( {x} \right)}}  \right]_{x = 0} =  \quad \left[
{\sum\limits_{k \ge 0} {\frac{{a_{k} b_{n - r}} }{{k_{\psi}  !\left( {n - r}
\right)_{\psi}  !}}n_{\psi}  !p_{0} \left( {x} \right)}}  \right]_{x = 0} =
\]
\[
=\sum\limits_{k \ge 0} {\left( {{\begin{array}{*{20}c}
 {n} \hfill \\
 {k} \hfill \\
\end{array}} } \right)_{\psi}  a_{k} b_{n - k}}
\]
\end{proof}

\begin{rem}{\em
In vain of the proof above one easily notices that the algebra $\Phi_{\psi}$
of $\psi $-exponential formal power series is isomorphic to the
reduced incidence algebra R(L(S)) where the isomorphism $\varphi $ is given
by the bijective correspondence:
\[
F_{\psi}  \left( {z} \right) = \sum\limits_{n \ge 0} {\frac{{a_{n}
}}{{n_{\psi}  !}}z^{n}\mathop { \to} \limits^{\varphi} }  f = \left\{
{f\left( {A,B} \right) = \left\{ {{\begin{array}{*{20}c}
 {a_{\left| {B - A} \right|} ;\;A \le B} \\
 {0;\quad otherwise} \\
\end{array}} ;\quad A,B \in L\left( {S} \right)} \right.} \right\}
\]
\noindent
where $f,g,h \in R(L(S))$ and $h:=f*g$ corresponds to $\psi
$-binomial convolution i.e. for
$$H_{\psi}  \left( {z} \right) = \sum\limits_{n \ge 0} {\frac{{c_{n}
}}{{n_{\psi}  !}}} \textrm{ and } G_{\psi}  \left( {z} \right) = \sum\limits_
{n \ge 0}{\frac{{b_{n}} }{{n_{\psi}  !}}z^{n}}\textrm{ ,  }[n]W(z) \equiv
 [n](f^*g)(z) \equiv c_{n}\, \textrm{we get}$$
\[
c_{n} = \sum\limits_{k \ge 0}^{n} {\left( {{\begin{array}{*{20}c}
 {n} \hfill \\
 {k} \hfill \\
\end{array}} } \right)_{\psi}  a_{k} b_{n - k}}  .
\]
This remark and the isomorphism $\Sigma _{\psi} \approx \Phi
_{\psi}  $ constitute expected link of finite operator $\psi $-calculus
with incident algebras . Namely- the $\psi $- extension of Rota`s operator
calculus is a general representation of the algebra structure of $R(L(S))$.\\
In any kind of algebra the principal question is to find out which elements
have their inverse. In the case of the algebra $\Sigma _{\psi}  $
we know the answer due to the isomorphism Theorem~\ref{ththreethree}.
Indeed, in the algebra $\Phi _{\psi}  $ of
$\psi $-exponential formal power series its` element
$\;\sum\limits_{k \ge 0} {\frac{{a_{k}} }{{k_{\psi
} !}}} x^{k}$ is invertible iff $a_{0} \ne 0$ because then and only then
the infinite system of linear equations $a_{0} b_{0} = 1\,;\,a_{0} b_{1} +
a_{1} b_{0} = 0\,;\,a_{0} b_{2} + a_{1} b_{1} + a_{2} b_{0} = 0$... has
solution and then this solution is unique.}
\end{rem}

\begin{cor}\label{corthreeone}
Operator $T_{\partial _{\psi} } \in \Sigma_{\psi}$ has its`
inverse $T_{\partial _{\psi} }  ^{ - 1}\; \in \Sigma_{\psi}$ iff
$T_{\partial_{\psi} }1 \ne 0$.
\end{cor}

\begin{proof}
Take $n=0$ in $\left[ {T_{\partial _{\psi} }  S_{\partial
_{\psi} }  p_{n} \left( {x} \right)} \right]_{x = 0} \;$= $\sum\limits_{k
\ge 0} {\left( {{\begin{array}{*{20}c}
 {n} \hfill \\
 {k} \hfill \\
\end{array}} } \right)_{\psi}  a_{k} b_{n - k}}  $ and use the Theorem%
~\ref{ththreethree}
where $T_{\partial _{\psi} }  = \sum\limits_{n \ge 0} {\frac{{a_{n}
}}{{n_{\psi}  !}}} \partial _{\psi} ^{n} $ and $S_{\partial _{\psi} }  =
\sum\limits_{n \ge 0} {\frac{{b_{n}} }{{n_{\psi}  !}}} \partial _{\psi} ^{n}
$.
\end{proof}

\begin{rem}
{\em The operator
$E^{a}\left( {\partial _{\psi} }  \right) = exp_{\psi}  \{ a\partial _{\psi
} \} $ is invertible in $\Sigma_{\psi}$ but it is not $\partial _{\psi
} $-delta operator. From Corollary~\ref{corthreeone} we infer that no one
of $\partial_{\psi}  $-delta operators $Q\left( {\partial _{\psi} }  \right)$
is invertible.\\
>From the stated above and the first expansion theorem we have the next
corollary.}
\end{rem}

\begin{cor}
Operator $R_{\partial _{\psi} } \in \Sigma_{\psi} $ is a $\partial
_{\psi}$-delta operator iff $a_{0} = 0$ and $a_{1} \ne 0$, where
$R_{\partial _{\psi} }  = \sum\limits_{n \ge 0} {\frac{{a_{n}} }{{n_{\psi}
!}}} Q\left( {\partial _{\psi} ^{}}  \right)^{n}$ or equivalently :
$r(0) = 0$ {\em \&} $r'(0) \ne 0$ where $r(x) =
\sum\limits_{k \ge 0} {\frac{{a_{k}}}{{k_{\psi}  !}}} x^{k}\;$
is the correspondent of $R_{\partial _{q}}\;$ under the isomorphism Theorem%
~\ref{ththreethree}.
\end{cor}

Now note that every $\partial _{\psi}  $-delta operator $Q\left( {\partial
_{\psi} } \right)$ is a function $Q$ of $\partial _{\psi}  $ with the
expansion
\[
Q\left( {\partial _{\psi} }  \right) = \sum\limits_{n \ge 1} {\frac{{q_{n}
}}{{n_{\psi}  !}}} \partial _{\psi} ^{n}
\]
while $\partial _{\psi} $-delta operator $\partial _{\psi}  $ is a function of
$\partial _{\psi} $-delta operator $Q\left( {\partial _{\psi} }  \right)$
with the expansion
\[
\partial _{\psi}  = Q^{ - 1}\left( {Q\left( {\partial _{\psi} }  \right)}
\right) = \sum\limits_{n \ge 0} {\frac{{q_{n}} }{{n_{\psi}  !}}} Q^{ -
1}\left( {Q\left( {\partial _{\psi} }  \right)} \right)^{n} \quad ;
\]
\begin{center}
$Q \circ Q^{-1} = id$.
\end{center}

We shall now derive an expression for the $\psi $-exponential generating
function for $\partial _{\psi}  $-basic polynomial sequence of the $\partial
_{\psi}  $-delta operator $Q\left( {\partial _{\psi} }  \right)$ using
notation established above.

\begin{rem}{\em
$exp_{\psi} \{zx\}$ is the $\psi$-exponential generating function
for $\partial _{\psi}  $-basic polynomial sequence $\left\{ {x^{n}}
\right\}_{n = 0}^{\infty}  $ of the $\partial _{\psi}  $ operator.}
\end{rem}

\begin{cor}
The $\psi $-exponential generating function for $\partial _{\psi}  $-basic
polynomial sequence $\left\{ {p_{n} \left( {x} \right)} \right\}_{n =
0}^{\infty}  $ of the $\partial _{\psi}  $-delta operator $Q\left( {\partial
_{\psi} }  \right)$ is given by the following formula
\[
\sum\limits_{k \ge 0} {\frac{{p_{k} \left( {x} \right)}}{{k_{\psi}  !}}}
z^{k}\; = exp_{\psi}  \{ xQ^{ - 1}\left( {z} \right)\}.
\]
\end{cor}

\begin{proof}
>From $E^{a}\left( {\partial _{\psi} }  \right) = \sum\limits_{k \ge 0}
{\frac{{a_{k}} }{{k_{\psi}  !}}} Q\left( {\partial _{\psi} }  \right)^{k}$
where $a_{k} = \left[ {E^{a}\left( {\partial _{\psi} }  \right)p_{k} \left(
{x} \right)} \right]_{x = 0} = p_{k} \left( {a} \right)$ we get
\[
E^{a}\left( {\partial _{\psi} }  \right) = \sum\limits_{k \ge 0}
{\frac{{p_{k} \left( {a} \right)}}{{k_{\psi}  !}}} Q\left( {\partial _{\psi
}}  \right)^{k}.
\]
Applying now the isomorphism Theorem~\ref{ththreethree} we have
\[
exp_{\psi}  \left( {xa} \right) = \sum\limits_{k \ge 0} {\frac{{p_{k} \left(
{a} \right)}}{{k_{\psi}  !}}} Q\left( {x} \right)^{k}.
\]
After the substitution $Q(x) = z$ and $a=x$ we get the thesis.
\end{proof}


\section{The Pincherle $\psi$-Derivative and Sheffer $\psi$-polynomials}

We again deliberately follow the way of Rota`s systematic presentation as the
finite operator calculus is already well suited for its
generalizations encompassing $\psi $-extension. We mostly use Rota`s
notation with indispensable changes.

\subsection{Pincherle $\psi $-derivative}

In order to define the Pincherle $\psi $-derivative ``\textbf{'}'' as a
linear map\\ \textbf{'} : $\Sigma_{\psi}  \quad  \to  \quad \Sigma_{\psi} $
we define $\hat {x}_{\psi}$.

\begin{defn} {\em (compare with (2) in \cite{2} and see also \cite{29} and
\cite{9b})}\\
$\hat {x}_{\psi}  :P \to P;
\quad
\hat {x}_{\psi}  x^{n} = \frac{{\psi _{n + 1} \left( {q} \right)\left( {n +
1} \right)}}{{\psi _{n} \left( {q} \right)}}x^{n + 1} =
\frac{{\left( {n + 1} \right)}}{{\left( {n + 1} \right)_{\psi} } }x^{n +
1},\quad n \ge 0$.
\end{defn}

\begin{defn} {\em (compare with (17) in \cite{2} )}\\
The Pincherle $\psi $-derivative is the linear map {\bf '} : $\Sigma _{\psi}
\quad \to  \quad \Sigma _{\psi} $ defined according to

$T_{\partial _{\psi} }$ \textbf{ '} $ = T_{\partial _{\psi} }  \;\hat
{x}_{\psi} - \hat {x}_{\psi}  T_{\partial _{\psi} }  \;
\equiv \quad $\textbf{[}$T_{\partial _{\psi} }  \;$\textbf{ ,} $\hat {x}
_{\psi}$\textbf{].}
\end{defn}

It is easy to see that $T_{\partial _{\psi} }$ {\bf '} $\in$
\textbf{\textit{ $\Sigma_{\psi} $}} for $T_{\partial _{\psi} }  \; \in
 \quad \Sigma_{\psi} $ i.e.\\$\forall $ a $ \in $\textit{F} ; [$E^{a}\left(
{\partial _{q}}  \right)$, [$T_{\partial _{\psi} }$, $\hat {x}_{\psi}
$]]$ = 0$
The later follows inmediately from \\
$[\partial _{\psi}  $, [$T_{\partial _{\psi}
} \;$ , $\hat {x}_{\psi}  $]]$_{} = 0$
and the expansion $E^{a}\left( {\partial _{\psi} }
\right) = \sum\limits_{k \ge 0} {\frac{{a^{k}}}{{k_{\psi}  !}}} \partial
_{\psi}^{k}$ due to $\partial _{\psi}$ {\bf '}$=1$
(as direct verification in the $\left\{ {x^{n}} \right\}_{n = 0}^{\infty}
$ basis shows). Therefore the linear map
{\bf '}: $\Sigma_{\psi}  \quad  \to  \quad \Sigma_{\psi} ${\bf '} is well
defined. Because $\partial _{\psi} ${\bf '}=1 one may introduce in a formal
sense the appealing notation \textbf{'} $
\equiv \frac{{d}}{{d\partial_{\psi}}}$. It is easy to see that the usual
rules hold for Pincherle $\psi $-derivative as for example
\begin{center}
($\partial _{\psi}  ^{n}$){\bf '} $= \partial _{\psi}  ^{n - 1} \quad
\Leftrightarrow \quad \frac{{d}}{{d\partial _{\psi} } }\partial _{\psi}
^{n} = n\partial_{\psi}  ^{n - 1}$.
\end{center}
Pincherle $\psi $-derivative is therefore the true derivation with Leibnitz
rule (see: Proposition \ref{propfourone})

\begin{defn} In accordance with the isomorphism theorem every $\psi$-shift
invariant operator $T_{\partial _{\psi} }  \in \Sigma _{\psi} $ of the form
$T_{\partial _{\psi} }  = \sum\limits_{n \ge 0} {\frac{{a_{n}} }{{n_{\psi}
!}}} \partial _{\psi} ^{n} $, has as its unique correspondent the $\psi
$-exponential formal power series
$t\left( {z} \right) = \sum\limits_{k \ge 0} {\frac{{a_{k}} }{{k_{\psi}  !}}}
z^{k}.$ We shall call $t\left( {z} \right)$ the \textit{indicator} of
$T_{\partial _{\psi}}$ operator.
\end{defn}

\begin{cor}\label{corfourone}
Let $t\left( {z} \right) = \sum\limits_{k \ge 0} {\frac{{a_{k}} }{{k_{\psi}
!}}} z^{k}$ be the indicator of $T_{\partial _{\psi} }  \in \Sigma _{\psi} $.
Then\\ $t'\left( {z} \right) = \sum\limits_{k \ge 1} {\frac{{ka_{k}
}}{{k_{\psi}  !}}} z^{k - 1}$
is the indicator of $T_{\partial _{\psi} }  \;${\bf '}$ \in \Sigma _{\psi}$.
\end{cor}

Due to the isomorphism theorem and the Corollary~\ref{corfourone}
the Leibniz rule holds which is obvious after one notices that
''\textbf{'}'' is the derivation in the commutative algebra
$\Sigma_{\psi} $ of $\partial _{\psi}  $-shift invariant linear
operators on the commutative algebra $P$ of polynomials.

\begin{prop}\label{propfourone}
($T_{\partial _{\psi} } \;S_{\partial _{\psi} }  \;$)\textbf{'} $=
T_{\partial _{\psi} } \;$\textbf{'} $S_{\partial _{\psi} }\; +
S_{\partial _{\psi} } \;T_{\partial _{\psi} }\;$\textbf{'}$\;$ ;
$T_{\partial _{\psi} } \;$, $S_{\partial_{\psi} } \; \in \; \Sigma _{\psi} $.
\end{prop}

As an immediate consequence of the Proposition~\ref{propfourone} we get
\begin{center}
($S_{\partial _{\psi} }  ^{n}\;$)\textbf{'}= n $S_{\partial _{\psi} }
\;$\textbf{'}$S_{\partial _{\psi} }  ^{n - 1}\; \quad \forall S_{\partial
_{\psi} }  \; \in  \quad \Sigma _\psi$.
\end{center}
Like in \cite{12} one infers from the isomorphism theorem that the following
is true.

\begin{prop}
$Q\left( {\partial _{\psi} }  \right)$ is the $\partial _{\psi}
$-delta operator iff there exists invertible\\
$S_{\partial_{\psi}} \in \Sigma_{\psi}$ such that
$$
Q\left( {\partial _{\psi} }  \right) \; = \;
\partial _{\psi}
S_{\partial _{\psi} }.
$$
\end{prop}
\begin{proof} $\textrm{  }$\\
$\Rightarrow $Let $S_{\partial _{\psi} } \; = \; \sum\limits_{k \ge 0}
{\frac{{s_{k}} }{{k_{\psi}  !}}} \partial _{\psi}  ^{k}$ and
$\exists  \quad S_{\partial_{\psi}}^{ - 1}$.
Then $\partial_{\psi}  S_{\partial_{\psi} }$ is a $\partial_{\psi}  $-delta
operator because $s_{0} \ne 0$.\\
\\$ \Leftarrow $ Let $Q\left( {\partial _{\psi} }  \right) = \sum\limits_{k
\ge 1} {\frac{{q_{k}} }{{k_{\psi}  !}}} \partial _{\psi}  ^{k}\quad, q_{1}
\ne 0$ be a $\partial _{\psi}  $-delta operator. Then \{
``\textit{S =Q/}$\partial _{\psi}  $'' \}
\[
S_{\partial _{\psi} } \; = \; \sum\limits_{k \ge 0} {\frac{{q_{k + 1}} }{{\left(
{k + 1} \right)_{\psi}  !}}} \partial _{\psi}  ^{k}
 \equiv
\quad
\sum\limits_{k \ge 0} {\frac{{s_{k}} }{{k_{\psi}  !}}} \partial _{\psi}
^{k}\quad ;\quad s_{0} = q_{1} \ne 0.
\]
\end{proof}

\begin{ex}{\em
$\Delta _{q} = \partial _{q} $\textit{S$\; \equiv \;
$}$\partial _{q} \sum\limits_{k \ge 0} {\frac{{\partial _{q} ^{k}}}{{\left(
{k + 1} \right)_{q} !}}}  = E^{a}\left( {\partial _{q}}  \right) - id$.\\

$\Delta _{\psi}  = \partial _{\psi}  $\textit{S$ \; \equiv \;$}$\partial
_{\psi} \sum\limits_{k \ge 0} {\frac{{\partial _{\psi}^{k}}}{{\left( {k + 1}
\right)_{\psi}  !}}}\; = \;E^{a}\left( {\partial _{\psi} }  \right) - id$.}
\end{ex}

The Pincherle $\psi $-derivative notion appears very effective in
formulating expressions for $\partial _{\psi}  $-basic polynomial sequences
of the given $\partial _{\psi}  $-delta operator $Q\left( {\partial _{\psi}
} \right)$. This is illustrated by (compare with 1.1.37. in \cite{9b} ) the
theorem that follows.

\begin{thm} \label{thfourone}
{\em ($\psi $-Lagrange and $\psi $-Rodrigues formulas)}\\
Let $\left\{ {p_{n} \left( {x} \right)} \right\}_{n = 0}^{\infty}  $ be
$\partial _{\psi}  $-basic polynomial sequence of the $\partial _{\psi}
$-delta operator $Q\left( {\partial _{\psi} }  \right)$:

 $Q\left( {\partial _{\psi} }  \right) \quad = \partial _{\psi}
S_{\partial _{\psi} }.$ Then for $n>0$:
\begin{enumerate}
\renewcommand{\labelenumi}{\em (\arabic{enumi})}
\item\label{one} $p_{n}(x) = Q\left( {\partial _{\psi} }
\right)$\textbf{'} $S_{\partial _{\psi} }^{-n-1}\;${\em x}$^{n}$ ;

\item\label{two} $p_{n}(x) = S_{\partial_{\psi} }^{-n}${\em x}$
^{n} - \frac{{n_{\psi} } }{{n}}$ ($S_{\partial _{\psi} }
^{ - n}\;$)\textbf{'}{\em x}$^{n-1};$

\item\label{three} $p_{n}(x) = \frac{{n_{\psi} } }{{n}}\hat {x}_{\psi}
S_{\partial _{\psi} }^{ - n}${\em x}$^{n-1}$;

\item\label{four} $p_{n}(x) = \frac{{n_{\psi} } }{{n}}\hat {x}_{\psi}
(Q\left( {\partial _{\psi} }  \right)$\textbf{'} )$^{-1}
p_{n-1}(x)$  {\em (Rodrigues $\psi $-formula )}.
\end{enumerate}
\end{thm}
\begin{proof}

Temporarily we use in the proof the following abbreviations: $Q\left(
{\partial _{\psi} }  \right) = Q$; $S_{\partial _{\psi} }\; = S$

\begin{enumerate}
\renewcommand{\labelenumi}{\bf \Roman{enumi}.}
\item We shall prove that right-hand sides of (\ref{one}) and (\ref{two})
determine the same polynomial sequence:

$Q$\textbf{'}$S^{-n-1} = (\partial _{\psi}S)$\textbf{'} $S^{-n-1} = ( S +
\partial _{\psi}S$\textbf{'} $)S^{-n-1} =
S^{-n} +S$\textbf{'}$S^{-n-1}\partial _{\psi}
= S^{-n} - (1/n) (S^{-n}$ )\textbf{'} $\partial _{\psi} $\\
due to this $Q$\textbf{'}$S^{-n-1}$ x$^{n}
= S^{-n}$ x$^{n} - \frac{{n_{\psi} } }{{n}}
(S^{-n}$ )\textbf{'}x$^{n-1}$

\item We now prove that right-hand sides of (\ref{one}) and (\ref{two}) and
(\ref{three}) determine the same polynomial sequence:

$S^{-n}$x$^{n} - \frac{{n_{\psi} } }{{n}}
(S^{-n}$)\textbf{'}x$^{n-1} =S ^{-n}$x$^{n} -
\frac{{n_{\psi} } }{{n}}(S^{-n} \hat {x}_{\psi} - \hat
{x}_{\psi} S^{-n}$)x$^{n-1} = \frac{{n_{\psi}
}}{{n}}\hat {x}_{\psi} S^{-n}$ x$^{n-1}$

\item We denote this polynomial sequence by
$q_{n}(x)= Q$\textbf{'}$S^{-n-1}$ x$^{n}$

We prove that the sequence $q_{n}(x)$ satisfies the requirements of the
$\partial _{\psi}  $-basic polynomial sequence

\begin{enumerate}
\renewcommand{\labelenumi}{(\alph{enumi})}
\item $Q q_{n}(x) = \partial _{\psi} SQ$\textbf{'}$S^{-n-1}$x$^{n} =
Q$\textbf{'}$S^{-n}\partial _{\psi}  $ x$^{n} = n_{\psi}q_{n-1}(x)$;

\item $q_{n}(x) = \frac{{n_{\psi} } }{{n}}\hat {x}_{\psi} S^{-n}$x$^{n-1}
\quad  \Rightarrow  \quad q_{n}(0)=0$; $n>0$;

\item $q_{n}(x) = \frac{{n_{\psi} } }{{n}}\hat {x}_{\psi
} S^{-n}$x$^{n-1} \quad  \Rightarrow \quad q_{0}$(x)$= 1$.
\end{enumerate}

Therefore in our notation : $q_{n}(x) \equiv p_{n}(x)$.

\item We now derive the Rodrigues formula. According to
Corollary~\ref{corfourone} and Corollary~\ref{threeone} there exists
(\textit{Q} ')$^{-1}$. Hence \textit{Q}'\textit{S}$^{-n-1}$ x $^{n}
= p_{n}$\textit{(x)} may be rewritten as x$^{n} =
S^{n+1}$ (\textit{Q} ')$^{-1}$ \textit{p}$_{n}$\textit{(x).}
Inserting x$^{n-1} = S^{n}$(\textit{Q} ')$^{-1}$
\textit{p}$_{n-1}$\textit{(x)} into (\ref{three}) $p_{n}(x) =
\frac{{n_{\psi} } }{{n}}\hat {x}_{\psi}  $\textit{} $S_{\partial _{q}}
^{ - n}\;$\textit{} x$^{n-1} $ one gets the Rodrigues formula (\ref{four})
$p_{n}(x) = \frac{{n_{\psi} } }{{n}}\hat {x}_{\psi}  $
( \textit{Q} ' )$^{-1}$ \textit{p}$_{n-1}$\textit{(x}).
\end{enumerate}
\end{proof}

\begin{com} \label{com1} {\em
Let $\left\{ {p_{n}}  \right\}_{n \ge 0} $ be the normal
polynomial sequence i.e. $p_{0} \left( {x} \right)$= 1, deg $p_{n}
= n$ and $p_{n} \left( {0} \right)$= 0 ; $n \ge 1$ \cite{29}. Then
we call $\left\{ {p_{n}}  \right\}_{n \ge 0} $ the $\psi $-basic
sequence of the operator $Q$ and $Q$ in its turn is called
generalized differential operator \cite{29} if in addition
$Q\,p_{n} = n_{\psi} p_{n - 1} $. Paralelly we define a linear map
$\hat {x}_{Q} $: $P \to P$ such that $\hat {x}_{Q} p_{n}
=\frac{{\left( {n + 1} \right)}}{{\left( {n + 1} \right)_{\psi} }
}p_{n + 1};\;n \ge 0$. One notices immediately that the Theorem
\ref{thfourone} holds also for the above operators $Q$ and $\hat
{x}_{Q} $ (see Theorem 4.3. in \cite{29} ). }
\end{com}

\begin{ex}{\em
Put in Rodrigues formula $Q\left( {\partial _{\psi} }  \right) = \partial
_{\psi}  $. Then \textit{p}$_{n}$\textit{(x) = x}$^{n},\\
n=0,1,2,...$ because $\partial _{\psi}  $\textbf{'} $= 1$ and $p_{0}(x) =1$.}
\end{ex}

\begin{cor}\label{corfourtwo}
Let $Q\left( {\partial _{\psi} } \right) = \partial _{\psi}
S_{\partial _{\psi} }$ and $R\left( {\partial _{\psi} }
\right) = \partial _{\psi}  P_{\partial _{\psi} } $ be
$\partial _{\psi}$-delta operators with $\partial _{\psi}  $-basic
polynomial sequences  $\left\{ {p_{n} \left( {x} \right)}
\right\}_{n = 0}^{\infty}  $ and $\left\{ {r_{n} \left( {x} \right)}
\right\}_{n = 0}^{\infty}  $ respectively. Then
\begin{enumerate}
\renewcommand{\labelenumi}{\em (\arabic{enumi})}
\setcounter{enumi}{4}
\item\label{five} $p_{n}(x) = R\left( {\partial _{\psi} }
\right)$\textbf{'}($Q\left( {\partial _{q}}  \right)$\textbf{'}
)$^{-1}S_{\partial _{\psi} }  ^{ - n - 1}\;P_{\partial _{\psi} }  ^{n +
1}r_{n}(x)$ $n \ge 0$;

\item \label{six}$p_{n}(x) = \hat {x}_{\psi} ( P_{\partial _{\psi}}\;
S_{\partial _{\psi} }^{ - 1})^{n}  \hat {x}_{\psi}^{ -1}r_{n}(x); n>0$.
\end{enumerate}
\end{cor}

\begin{proof}
ad.(\ref{five}) use (\ref{one}) from the Theorem~\ref{thfourone} for obvious
substitution; ad.(\ref{six}) use (\ref{three}) from the
Theorem~\ref{thfourone} for obvious substitution.
\end{proof}
To this end we shall present another characterization of $\partial _{\psi}
$-basic polynomial sequences of $\partial _{\psi}  $-delta operators
$Q\left( {\partial _{\psi} }  \right)$.

\begin{thm}
Let $P_{\partial _{\psi} } \in \Sigma_{\psi}$ be an invertible
$\partial _{\psi}  $-shift invariant operator. Let $\left\{ {p_{n} \left(
{x} \right)} \right\}_{n = 0}^{\infty}  $ be $\partial _{\psi}  $-basic
polynomial sequence of the $\partial _{\psi}  $-delta operator $Q\left(
{\partial _{\psi} }  \right)$ satisfying
$$
\left[ {x^{ - 1}p_{n} \left( {x} \right)} \right]_{x = 0} \;=
\;
n_{\psi}
\left[ {P_{\partial _{q}}  ^{ - 1}p_{n - 1} \left( {x} \right)} \right]_{x =
0} \; n>0.$$

Then $\left\{ {p_{n} \left( {x} \right)} \right\}_{n = 0}^{\infty}  $ is
$\partial _{\psi}  $-basic polynomial sequence of the $\partial _{\psi}
$-delta operator of the form: $Q\left( {\partial _{\psi} }
\right) = \partial _{\psi}  P_{\partial _{\psi} }\;$.
\end{thm}

\begin{proof}
The linear operator $Q\left( {\partial _{\psi} }  \right)$ is defined
completely by $Q\left( {\partial _{\psi} }  \right)1 = 0$ and
$Q\left( {\partial _{\psi} }  \right) p_{n}(x) = n_{\psi}p_{n-1}(x)$.
It is easy to see that $Q\left(
\partial _{\psi} \right) \in \Sigma_{\psi}$.
The condition
$\left[ {x^{ - 1}p_{n} \left( {x} \right)} \right]_{x = 0}
 = n_{\psi}  \left[ {P_{\partial _{\psi} }  ^{ - 1}p_{n - 1} \left( {x}
\right)} \right]_{x = 0} \;$ $n>0$ may be rewritten as follows:
 $\left[ {x^{ - 1}p_{n} \left( {x} \right)} \right]_{x = 0} \;= \quad
\left[ {P_{\partial _{\psi} }  ^{ - 1}Q\left( {\partial _{\psi} }
\right)p_{n} \left( {x} \right)} \right]_{x = 0} \;$ $n>0$.
By linearity argument the above extends to arbitrary polynomial i.e. we
have\\
 $\left[ {x^{ - 1}p\left( {x} \right)} \right]_{x = 0} \;=
\quad
\left[ {P_{\partial _{\psi} }  ^{ - 1}Q\left( {\partial _{\psi} }
\right)p\left( {x} \right)} \right]_{x = 0} \;$ $n>0$ .
At the same time for all polynomial such that $p(0) = 0$ the obvious
identity holds:
\[
\left[ {x^{ - 1}p\left( {x} \right)} \right]_{x = 0} \;=
\quad
\left[ {\partial _{\psi}  p\left( {x} \right)} \right]_{x = 0} \;
\]
due to which we have the identity
\[
\left[ {P_{\partial _{\psi} }  ^{ - 1}Q\left( {\partial _{\psi} }
\right)p\left( {x} \right)} \right]_{x = 0} \;=
\quad
\left[ {\partial _{\psi}  p\left( {x} \right)} \right]_{x = 0} \;
\]
satisfied now for all polynomials (including those with $p(0) \ne 0$ )\\
Put now $p(x) = q$(x $+_{\psi}a)$ and use the fact :
$P_{\partial _{\psi} }  \;$\textit{ ,} $Q\left( {\partial _{\psi} }
\right) \quad  \in \Sigma_{\psi} $ in order to obtain
\begin{multline*}
 \partial _{\psi}q(a) = \left[ {P_{\partial _{\psi} }^{ -
1}Q\left( {\partial _{\psi} }  \right)E^{a}\left( {\partial _{\psi} }
\right)q\left( {x} \right)} \right]_{x = 0} \; = \quad \left[ {E^{a}\left(
{\partial _{\psi} }  \right)P_{\partial _{\psi} }  ^{ - 1}Q\left( {\partial
_{\psi} }  \right)q\left( {x} \right)} \right]_{x = 0} \; = \\
=P_{\partial_{\psi} }  ^{ - 1}\;Q\left( {\partial _{\psi} }  \right)q(a)
\end{multline*}
$\forall $a$ \in F$. Hence $Q\left( {\partial _{\psi} }
\right) = \partial _{\psi}  P_{\partial _{\psi} }$.
\end{proof}

\begin{conc}{\em
One already sees that Rota`s finite operator calculus is ready and
well suited for natural almost automatic generalizations
encompassing $\psi $-extensions (``$\psi$-deformations'') and
$Q$-extensions (see Comment \ref{com1}). We have seen that the use
Rota`s calculus notion and and ideas of proofs with indispensable
changes leads to this kind of extension.

We confirm it further considering a very important notion of Sheffer $\psi
$-po\-ly\-no\-mials or correspondingly generalized Sheffer $\psi$-polynomials
- see Comment \ref{com1} and \cite{29}.

This later notion provides ultimate generalization of umbral calculus
\cite{7,10} in that sense that we deal now with all $\left\{ {p_{n}}  \right\}
_{o}^{\infty} $, deg $p_{n}=n$ polynomial sequences (see: \cite{29} and
\cite{38} ). Umbral calculus for other functions (not necessarily polynomials)
is well known and has vast literature (see: recent one \cite{38})}
\end{conc}

\subsection{Sheffer $\psi $-polynomials}

\begin{defn}
A polynomial sequence $\left\{ {s_{n} \left( {x} \right)} \right\}_{n =
0}^{\infty}  $ is called the Sheffer $\psi $-po\-ly\-no\-mials sequence of
the $\partial _{\psi}  $-delta operator $Q\left( {\partial _{\psi} }  \right)$
iff
\begin{enumerate}
\renewcommand{\labelenumi}{\em (\arabic{enumi})}
\item $s_{0} \left( {x} \right) = c \ne 0$,
\item $Q\left( {\partial_{\psi} }  \right)s_{n} \left( {x} \right) =
n_{\psi}  s_{n -1} \left( {x} \right)$.
\end{enumerate}
\end{defn}

The following proposition relates Sheffer $\psi $-polynomials of the
$\partial _{\psi}  $-delta operator $Q\left( {\partial _{\psi} }  \right)$
to the unique $\partial _{\psi}  $-basic polynomial sequence of $Q\left(
{\partial _{\psi} }  \right)$.

\begin{prop}\label{propfourthree}
Let $Q\left( {\partial _{\psi} }  \right)$ be a $\partial _{\psi}  $-delta
operator with $\partial _{\psi}  $-basic polynomial sequence $\left\{ {q_{n}
\left( {x} \right)} \right\}_{n = 0}^{\infty}  $. Then $\left\{ {s_{n}
\left( {x} \right)} \right\}_{n = 0}^{\infty}  $ is a sequence of Sheffer
$\psi $-polynomials of $Q\left( {\partial _{\psi} }  \right)$ iff there
exists an invertible $\partial _{\psi}  $-shift invariant operator
$S_{\partial _{\psi} }$ such that $s_{n} \left( {x}
\right)$=$S_{\partial _{\psi} }  ^{ - 1}\;q_{n} \left( {x} \right)$.
\end{prop}

\begin{proof}$\textrm{  }$\\
$ \Rightarrow \;$ Let there exist a $\partial _{\psi}  $-shift invariant
operator $S_{\partial _{\psi} }$ such that $s_{n} \left( {x}
\right) = S_{\partial _{\psi} }  ^{ - 1}\;q_{n} \left( {x} \right)$.
Of course [$S_{\partial _{\psi} }  ^{ - 1}\;$,$Q\left( {\partial _{\psi} }
\right)$] $= 0$. Hence $Q\left( {\partial _{\psi} }  \right)s_{n} \left(
{x} \right)$= $n_{\psi}  s_{n - 1} \left( {x} \right)$. Indeed;\\ $Q\left(
{\partial _{\psi} }  \right)s_{n} \left( {x} \right) = Q\left(
{\partial _{\psi} }  \right)S_{\partial _{\psi} }  ^{ - 1}\;q_{n} \left(
{x} \right) = S_{\partial _{\psi} }  ^{ - 1}\;Q\left( {\partial
_{\psi} }  \right)q_{n} \left( {x} \right) = n_{\psi}  s_{n -
1} \left( {x} \right)$. Also $s_{0} \left( {x} \right) = S_{\partial
_{\psi} }  ^{ - 1}\;q_{n} \left( {x} \right) = S_{\partial
_{\psi} }  ^{ - 1}\;1 = $c $ \ne 0$.\\
$ \Leftarrow \;$ Let $\left\{ {s_{n} \left( {x} \right)} \right\}_{n =
0}^{\infty}  $ be a Sheffer $\psi $- sequence of $Q\left( {\partial _{\psi}
} \right)$. The linear operator $S_{\partial _{\psi} }$ is defined via:
$S_{\partial _{\psi} }\; : \;s_{n} \left( {x} \right) \to q_{n} \left( {x}
\right)$. $S_{\partial _{q}}$ is of course invertible because deg $s_{n}
\left( {x} \right) =$ deg $q_{n} \left( {x} \right)$ and $s_{0} \left( {x}
\right) = c \ne 0$. Also $S_{\partial _{\psi} }  \in \; \Sigma_{\psi}$.
Indeed. [$S_{\partial _{\psi} }  \;$, $Q\left( {\partial _{\psi} }
\right)$] $= 0$, as easily checked in $\left\{ {s_{n} \left( {x} \right)}
\right\}_{n = 0}^{\infty}  $ basis. Then [$E^{a}\left( {\partial
_{\psi} }  \right)$, $S_{\partial _{\psi} }  \;$] $= 0$ , $\forall a \in F$
because $E^{a}\left( {\partial _{\psi} }  \right) = \sum\limits_{n \ge 0} {\frac{{a_{n}} }{{n_{\psi}  !}}}
Q\left( {\partial _{\psi}}  \right)^{n}$.
\end{proof}

\begin{conc}{\em
The family of Sheffer $\psi $-polynomials`
sequences $\left\{ {s_{n} \left( {x} \right)} \right\}_{n = 0}^{\infty}  $
corresponding to a fixed $\partial _{\psi}  $-delta operator $Q\left(
{\partial _{\psi} }  \right)$ is labeled by the abelian group of all
invertible operators $S_{\partial _{\psi} } \in  \; \Sigma_{\psi} $.
This families of Sheffer $\psi$-polynomials are orbits of such groups
contained in the algebra $\Sigma_{\psi} $.\\
\textbf{Naming}: Sheffer $\psi $-polynomial sequence $\left\{ {s_{n} \left(
{x} \right)} \right\}_{n = 0}^{\infty}  $ labeled by $S_{\partial _{\psi} }
\;$ is refered to as the sequence of Sheffer $\psi $-polynomials of the
$\partial _{\psi}  $-delta operator $Q\left( {\partial _{\psi} }  \right)$
\textit{relative to} $S_{\partial _{\psi} }$.}
\end{conc}

The following theorem is valid for invertible $S_{\partial _{\psi} }$ such
that the sequence of Sheffer $\psi $-polynomials of the $\partial _{\psi}
$-delta operator $Q\left( {\partial _{\psi} }  \right)$ is relative to
$S_{\partial _{\psi} }$.

\begin{thm}\label{thfourthree} {\em (Second $\psi $- Expansion Theorem)}\\
Let $Q\left( {\partial _{\psi} }  \right)$ be the $\partial _{\psi}  $-delta
operator $Q\left( {\partial _{\psi} }  \right)$ with the $\partial _{\psi}
$-basic polynomial sequence $\left\{ {q_{n} \left( {x} \right)} \right\}_{n
= 0}^{\infty}  $. Let $S_{\partial _{\psi} }$ be an \textit{invertible}
$\partial _{\psi}  $-shift invariant operator and let $\left\{ {s_{n} \left(
{x} \right)} \right\}_{n = 0}^{\infty}  $ be its sequence of Sheffer $\psi
$-polynomials. Let $T_{\partial _{\psi} }  \;$be \textit{any} $\partial
_{\psi}  $-shift invariant operator and let \textit{p(x)} be any polynomial.
Then the following identity holds :

\begin{center}
$\forall \; y \in F \wedge \; \forall \; p \in P \; \quad T_{\partial
_{\psi} }  \;p\left( {x + _{\psi}  y} \right) = \sum\limits_{k \ge 0}
{\frac{{s_{k} \left( {y} \right)}}{{k_{\psi}  !}}} Q\left( {\partial
_{\psi} }  \right)^{k}S_{\partial _{\psi} }  \;T_{\partial _{\psi} }
\;p\left( {x} \right)$ .
\end{center}
\end{thm}

\begin{proof}
According to the first expansion theorem
\[
E^{y}\left( {\partial _{\psi} }  \right) \quad =
\quad
\sum\limits_{k \ge 0} {\frac{{q_{k} \left( {y} \right)}}{{k_{\psi}  !}}}
Q\left( {\partial _{\psi} }  \right)^{k}
\]
from which - by acting on $p\left( {x} \right)$ - we have

 $p\left( {x + _{\psi}  y} \right)=
\sum\limits_{k \ge 0} {\frac{{q_{k} \left( {y} \right)}}{{k_{\psi}  !}}}
Q\left( {\partial _{\psi} }  \right)^{k}
p\left( {x} \right)$ or equivalently\\

 $p\left( {x + _{\psi}  y}
\right) = \sum\limits_{k \ge 0} {\frac{{q_{k} \left( {x} \right)}}{{k_{\psi
} !}}} Q\left( {\partial _{\psi} }  \right)^{k}p\left( {y} \right)$ .

Of course $\quad S_{\partial _{\psi} }  ^{ - 1}\;p\left( {x + _{\psi}  y}
\right) = \sum\limits_{k \ge 0} {\frac{{s_{k} \left( {x} \right)}}{{k_{\psi
} !}}} Q\left( {\partial _{\psi} }  \right)^{k}p\left( {y} \right) \quad$ or
equivalently
\[
p\left( {x + _{\psi}  y} \right)=
\sum\limits_{k \ge 0} {\frac{{s_{k} \left( {y} \right)}}{{k_{\psi}  !}}}
Q\left( {\partial _{\psi} }  \right)^{k}
S_{\partial _{\psi} }  \;
p\left( {x} \right).
\]
Application of $T_{\partial _{\psi} }  \in \;\Sigma_{\psi} $ to
both sides yields the thesis.
\end{proof}

The theorem (\ref{thfourthree}) - called after Rota \cite{12} The Second
Expansion Theorem - has as an important outcome the following Corollary.

\begin{cor} {\em ($\psi $- expansion)}\\
Let $\left\{ {s_{n} \left( {x} \right)} \right\}_{n = 0}^{\infty}  $ be the
sequence of Sheffer $\psi $-polynomials of the $\partial _{\psi}  $-delta
operator $Q\left( {\partial _{\psi} }  \right)$ relative to $S_{\partial
_{\psi} }$. Then
\[
S_{\partial _{\psi} }  ^{ - 1}\;=
\quad
\sum\limits_{k \ge 0} {\frac{{s_{k} \left( {0} \right)}}{{k_{\psi}  !}}}
Q\left( {\partial _{\psi} }  \right)^{k}.
\]
\end{cor}

\begin{proof}
Set $y = 0$ and take $T_{\partial _{\psi} }  \; = S_{\partial
_{\psi} }  ^{ - 1}\;$in the second expansion theorem thesis.
\end{proof}

The property of $\psi $-binomiality of $\partial _{\psi}  $-basic polynomial
sequence of the $\partial _{\psi}  $-delta operator $Q\left( {\partial
_{\psi} }  \right)$ has its straightforward counterpart for $Q\left(
{\partial _{\psi} }  \right)$`s Sheffer $\psi $-sequence relative to
$S_{\partial _{\psi} }$.

\begin{thm}\label{thfourfour} {\em (The Sheffer $\psi $-Binomial Theorem)}\\
Let $Q\left( {\partial _{\psi} }  \right)$ be the $\partial _{\psi}  $-delta
operator with the $\partial _{\psi}  $-basic polynomial sequence $\left\{
{q_{n} \left( {x} \right)} \right\}_{n = 0}^{\infty}  $. Let $\left\{
{s_{n} \left( {x} \right)} \right\}_{n = 0}^{\infty}  $ be the sequence of
Sheffer $\psi $-polynomials of $Q\left( {\partial _{\psi} }  \right)$
relative to\textbf{} an invertible $\partial _{\psi}  $-shift invariant
operator $S_{\partial _{\psi} }$. Then the following identity is true
\[
s_{n} \left( {x + _{\psi}  y} \right)=
\sum\limits_{k \ge 0} {\left( {{\begin{array}{*{20}c}
 {n} \hfill \\
 {k} \hfill \\
\end{array}} } \right)} _{\psi}  s_{k} \left( {x} \right)q_{n - k} \left(
{y} \right).
\]
\end{thm}

\begin{proof}
In order to get the thesis apply $S_{\partial _{\psi} }  ^{ - 1}\;$ to
both sides of $\psi $-binomial formula
\[
q_{n} \left( {x + _{\psi}  y} \right)=
\sum\limits_{k \ge 0} {\left( {{\begin{array}{*{20}c}
 {n} \hfill \\
 {k} \hfill \\
\end{array}} } \right)} _{\psi}  q_{k} \left( {x} \right)q_{n - k} \left(
{y} \right).
\]
\end{proof}

Analogously to the undeformed case \cite{12} Sheffer $\psi $-polynomials are
completely determined by their constant terms as seen from the successive
Corollary.

\begin{cor}
For $\left\{ {q_{n} \left( {x} \right)} \right\}_{n
= 0}^{\infty}  $ and $\left\{ {s_{n} \left( {x} \right)} \right\}_{n =
0}^{\infty}  $ defined as in Theorem (\ref{thfourfour}) the following
holds $s_{n} \left({x} \right) =
\sum\limits_{k \ge 0} {\left( {{\begin{array}{*{20}c}
 {n} \hfill \\
 {k} \hfill \\
\end{array}} } \right)} _{\psi}  s_{k} \left( {0} \right)q_{n - k} \left(
{x} \right)$.
\end{cor}
\begin{proof} Obvious. \end{proof}

Analogously to the undeformed case \cite{12} also the converse of the second
$\psi $-expansion theorem is true.

\begin{prop}\label{propfourfour}
Let $Q\left( {\partial _{\psi} }  \right)$ be a $\partial _{\psi}  $-delta
operator. Let $S_{\partial _{\psi} }$ be an invertible $\partial _{\psi
} $-shift invariant operator. Let $\left\{ {s_{n} \left( {x} \right)}
\right\}_{n = 0}^{\infty}  $ be a polynomial sequence. Let
\begin{center}
$\forall \; a \in F \wedge \; \forall \; p \in P \quad E^{a}\left(
{\partial _{\psi} }  \right)p\left( {x} \right) = \sum\limits_{k \ge 0}
{\frac{{s_{k} \left( {a} \right)}}{{k_{\psi}  !}}} Q\left( {\partial
_{\psi} }  \right)^{k}S_{\partial _{\psi} }  \;p\left( {x} \right)$ .
\end{center}
Then the polynomial sequence $\left\{ {s_{n} \left( {x} \right)} \right\}_{n
= 0}^{\infty}  $ is the sequence of Sheffer $\psi $-polynomials of the
$\partial _{\psi}  $-delta operator $Q\left( {\partial _{\psi} }  \right)$
relative to $S_{\partial _{\psi} }$.
\end{prop}

\begin{proof}
The above formula may be recasted in the form
\[
E^{a}\left( {\partial _{\psi} }  \right)
p\left( {x} \right)=
\sum\limits_{k \ge 0} {\frac{{S_{\partial _{\psi} }  s_{k} \left( {a}
\right)}}{{k_{\psi}  !}}}
Q\left( {\partial _{\psi} }  \right)^{k}
p\left( {x} \right)
\]
whereupon we take now $p\left( {x} \right) = q_{n} \left( {x} \right)$ and
where $\left\{ {q_{n} \left( {x} \right)} \right\}_{n = 0}^{\infty}  $ is
the the $\partial _{\psi}  $-basic polynomial sequence of $Q\left( {\partial
_{\psi} }  \right)$ thus arriving at
\[
q_{n} \left( {x + _{\psi}  a} \right)=
\sum\limits_{k \ge 0} {\frac{{S_{\partial _{\psi} }  s_{k} \left( {a}
\right)}}{{k_{\psi}  !}}}
n_{\psi} ^{\underline {k}}
q_{n - k} \left( {x} \right)
\]
where we set $x=0$ thus getting
\[
q_{n} \left( {a} \right)=
\sum\limits_{k \ge 0} {\left( {{\begin{array}{*{20}c}
 {k} \hfill \\
 {n} \hfill \\
\end{array}} } \right)} _{\psi}
S_{\partial _{\psi} }  \;
s_{n} \left( {a} \right)
q_{n - k} \left( {0} \right)
\]
hence $\forall $ a$ \in F \quad q_{n} \left( {a} \right) =
S_{\partial _{\psi} }  \;s_{n} \left( {a} \right)$.
\end{proof}

Similarily as in nonextended case \cite{12} further correspondent constructs,
propositions and theorems hold also for the $\psi $- extension of the finite
operator calculus. For example due to the above Proposition
\ref{propfourfour} we get another proof of the Rodrigues formula
\begin{center}
\textit{p}$_{n}$\textit{(x) =} $\frac{{n_{\psi} } }{{n}}\hat {x}_{\psi}  $
($Q\left( {\partial _{\psi} }  \right)$\textbf{'} )$^{-1}$
\textit{p}$_{n-1}$\textit{(x).}
\end{center}

\begin{prop}
Let $\left\{ {p_{n} \left( {x} \right)} \right\}_{n = 0}^{\infty}  $ be the
$\partial _{\psi}  $-basic polynomial sequence of the $\partial _{\psi}
$-delta operator $Q\left( {\partial _{\psi} }  \right)$. Then
\begin{center}
$p_{n}(x) = \frac{{n_{\psi} } }{{n}}\hat {x}_{\psi}
(Q\left( {\partial _{\psi} }  \right)$\textbf{'} )$^{-1}
p_{n-1}(x)$.
\end{center}
\end{prop}

\begin{proof}
Due to $E^{a}\left( {\partial _{\psi} }  \right) = \sum\limits_{k \ge 0}
{\frac{{a_{k}} }{{k_{\psi}  !}}} Q\left( {\partial _{\psi} }  \right)^{k}$
where $a_{k} = \left[ {E^{a}\left( {\partial _{\psi} }  \right)p_{k} \left(
{x} \right)} \right]_{x = 0} = p_{k} \left( {a} \right)$ we have
\[
\sum\limits_{n \ge 0} {\frac{{a^{n}}}{{n_{\psi}  !}}} \partial _{\psi} ^{n}
= E^{a}\left( {\partial _{\psi} }  \right) = \sum\limits_{k \ge 0}
{\frac{{p_{k} \left( {a} \right)}}{{k_{\psi}  !}}} Q\left( {\partial _{\psi
}}  \right)^{k}.
\]
The Pincherle derivative beeing applied to both sides yields
\begin{multline*}
a\sum\limits_{k \ge 0} {\frac{{a^{k}\partial _{\psi} ^{k}} }{{k_{\psi}
!}}} \frac{{\left( {k + 1} \right)}}{{\left( {k + 1} \right)_{\psi}
}} =\hat {a}_{\psi}  \sum\limits_{k \ge 0}
{\frac{{a^{k}\partial _{\psi} ^{k}} }{{k_{\psi}  !}}} = \\ = \hat {a}_{\psi}
E^{a}\left( {\partial _{\psi} }  \right)=\sum\limits_{k \ge 0}
{\frac{{p_{k + 1} \left( {a} \right)}}{{k_{\psi}  !}}\frac{{\left( {k + 1}
\right)}}{{\left( {k + 1} \right)_{\psi} } }} Q\left( {\partial _{\psi} }
\right)^{k}Q\left( {\partial _{\psi} }  \right)\textrm{\bf {'}}.
\end{multline*}
Compare this with $E^{a}\left( {\partial _{\psi} }  \right) = \sum\limits_{k
\ge 0} {\frac{{s_{k} \left( {a} \right)}}{{k_{\psi}  !}}} Q\left(
{\partial _{\psi} }  \right)^{k}S_{\partial _{\psi} }$ (Proposition
\ref{propfourfour}) and recall that\\
 $\hat {x}_{\psi}  x^{n}=
\frac{{\left( {n + 1} \right)}}{{\left( {n + 1} \right)_{\psi} } }x^{n +
1};\quad n \ge 0$ and $\hat {x}_{\psi}  ^{ - 1}x^{n}=\frac{{n_{\psi}
}}{{n}}x^{n - 1}; \quad n > 0$ in order to notice that from $E^{a}\left(
{\partial _{\psi} }  \right) = \sum\limits_{k \ge 0} {\frac{{\hat {a}_{\psi
} ^{ - 1}p_{k + 1} \left( {a} \right)}}{{k_{\psi}  !}}\frac{{\left( {k + 1}
\right)}}{{\left( {k + 1} \right)_{\psi} } }} Q\left( {\partial _{\psi} }
\right)^{k}Q\left( {\partial _{\psi} }  \right)$\textbf{'}  and from the
Proposition 4.4. it follows that $\{ \frac{{\left( {k + 1} \right)}}{{\left(
{k + 1} \right)_{\psi} } }\hat {x}_{\psi}  ^{ - 1}p_{k + 1} \left( {x}
\right)\} _{k = 0}^{\infty}  $ is a sequence of Sheffer $\psi $-polynomials
of $\partial _{\psi}  $-delta operator $Q\left( {\partial _{\psi} }
\right)$ - a sequence relative to the invertible $\partial _{\psi}  $-shift
invariant operator $Q\left( {\partial _{\psi} }  \right)$\textbf{'}. Hence

 $\frac{{n}}{{n_{\psi} } }\hat {x}_{\psi}  ^{ - 1}$ \textit{p}$_{n}$\textit{(x)
=} ($Q\left( {\partial _{\psi} }  \right)$\textbf{'} )$^{-1}$
\textit{p}$_{n-1}$\textit{(x)}  according to the proposition
\ref{propfourthree}.
\end{proof}
In what follows let us remark that the notion and properties of umbral
composition, cross-sequences etc \cite{12} straihtfowardly spread also to the
$\psi $-extended calculus case.\\
Here - for the sake of future examples - we shall formulate few more
propositions for which - with indispensable modifications - proofs of Rota
can be adopted.

\begin{prop}
Let $Q\left( {\partial _{\psi} }  \right)$ be a $\partial _{\psi}  $-delta
operator. Let $S_{\partial _{\psi} }$ be an invertible $\partial _{\psi
} $-shift invariant operator. Let \textit{q(t)} and s\textit{(t)}  be the
indicators of $Q\left( {\partial _{\psi} }  \right)$ and $S_{\partial _{\psi
}}$ operators. Let \textit{q$^{-1}$(t )}
be the inverse $\psi $-exponential formal power series inverse to
\textit{q(t)}. Then the $\psi $-exponential generating function of Sheffer
$\psi $-polynomials sequence $\left\{ {s_{n} \left( {x} \right)} \right\}_{n
= 0}^{\infty}  $ of $Q\left( {\partial _{\psi} }  \right)$ relative
to $S_{\partial _{\psi} }  \;$is given by
\[
\;\sum\limits_{k \ge 0} {\frac{{s_{k} \left( {x}
\right)}}{{k_{\psi}  !}}} z^{k}\; = \;{\frac{1}{s\left( {q^{ -
1}\left( {z} \right)} \right)}}\;exp_{\psi}  \{ xq^{ - 1}\left(
{z} \right)\}.
\]
\end{prop}
In addition we have the following characterization of Sheffer $\psi
$-polynomials:

\begin{prop}
A sequence $\left\{ {s_{n} \left( {x} \right)} \right\}_{n = 0}^{\infty}  $
is the sequence of Sheffer $\psi $-polynomials of the $\partial _{\psi}
$-delta operator $Q\left( {\partial _{\psi} }  \right)$ with the $\partial
_{\psi}  $-basic polynomial sequence $\left\{ {p_{n} \left( {x} \right)}
\right\}_{n = 0}^{\infty}  $ iff
\[
s_{n} \left( {x + _{\psi}  y} \right)=
\sum\limits_{k \ge 0} {\left( {{\begin{array}{*{20}c}
 {n} \hfill \\
 {k} \hfill \\
\end{array}} } \right)} _{\psi}  s_{k} \left( {x} \right)q_{n - k} \left(
{y} \right).
\]
\end{prop}
Also the recurrence formula from \cite{12} $\psi $-extends straightforwardly.

\begin{prop}
Let $\left\{ {p_{n} \left( {x} \right)} \right\}_{n = 0}^{\infty}  $ be a
polynomial sequence with $p_{0} \left( {x} \right) = c \ne 0$. If in adition
it is Sheffer $\psi $-sequence then for any $\partial _{\psi}  $-delta
operator A there exists the sequence of \textit{constants} $\left\{ {s_{n}}
\right\}_{n = 0}^{\infty}  $such that
\begin{center}
A$p_{n} \left( {x} \right)$ = $\sum\limits_{k \ge 0} {\left(
{{\begin{array}{*{20}c}
 {n} \hfill \\
 {k} \hfill \\
\end{array}} } \right)} _{\psi}  p_{k} \left( {x} \right)s_{n - k} \quad n
\ge 0$.
\end{center}
If the above recurrence holds for some $\partial _{\psi}  $-delta operator A
and some sequence of \textit{constants} $\left\{ {s_{n}}  \right\}_{n =
0}^{\infty}  $ then $\left\{ {p_{n} \left( {x} \right)} \right\}_{n =
0}^{\infty}  $ is a Sheffer $\psi $-sequence. (Operator A is not necessarily
associated with $\left\{ {p_{n} \left( {x} \right)} \right\}_{n = 0}^{\infty
} $).
\end{prop}

\begin{proof}
Adopt Rota`s proof \cite{12} - with obvious replacements - and use the identity
\[
\left( {n - k} \right)_{\psi}  \left( {{\begin{array}{*{20}c}
 {n} \hfill \\
 {k} \hfill \\
\end{array}} } \right)_{\psi}  = n_{\psi}  \left( {{\begin{array}{*{20}c}
 {n - 1} \hfill \\
 {\; k} \hfill \\
\end{array}} } \right)_{\psi}.
\]
\end{proof}
In effect we have another characterization of Sheffer $\psi $-polynomials.

\begin{cor}
A sequence $\left\{ {p_{n} \left( {x} \right)} \right\}_{n = 0}^{\infty}  $
is the sequence of Sheffer q-polynomials of the $\partial _{\psi}  $-delta
operator $Q\left( {\partial _{\psi} }  \right)$ with the $\partial _{\psi}
$-basic polynomial sequence $\left\{ {q_{n} \left( {x} \right)} \right\}_{n
= 0}^{\infty}  $ iff there exists such a $\partial _{\psi}  $-delta operator
A (not necessarily associated with $\left\{ {p_{n} \left( {x} \right)}
\right\}_{n = 0}^{\infty}  $) and the sequence $\left\{ {s_{n}}  \right\}_{n
= 0}^{\infty}$ of \textit{constants} such that
\begin{center}
A$p_{n} \left( {x} \right)$ = $\sum\limits_{k \ge 0} {\left(
{{\begin{array}{*{20}c}
 {n} \hfill \\
 {k} \hfill \\
\end{array}} } \right)} _{\psi}  p_{k} \left( {x} \right)s_{n - k}, \quad n
\ge 0$.
\end{center}
\end{cor}
As in the undeformed case the natural inner product may be assosiated with
the sequence $\left\{ {s_{n} \left( {x} \right)} \right\}_{n = 0}^{\infty}
$ of Sheffer $\psi $-polynomials of the $\partial _{\psi}  $-delta operator
$Q\left( {\partial _{\psi} }  \right)$ relative to\textbf{} an invertible
$\partial _{\psi}  $-shift invariant operator $S_{\partial _{\psi} }  \;$.
For that purpose define the linear umbral operator
W : $s_{n} \left( {x} \right) \to \; x^{n}$.

\begin{defn}
Let $S_{\partial _{\psi} }$ be a $\partial _{\psi}  $-shift invariant
operator. Let W be the linear operator such that

W : $s_{n} \left( {x} \right) \to  \; x^{n}$ and then extended by
linearity. We define the following bilinear form
\begin{center}
(f(x),g(x))$\psi $ := \textbf{[(}Wf\textbf{)}\textbf{(} $Q\left( {\partial
_{\psi} } \right) \textbf{)}S_{\partial _{\psi} }  \;$g(x)\textbf{]}$_{x=0}$ f,g $ \in P$.
\end{center}
\end{defn}
It is easy to state the important property of this bilinear form now on
the reals.

\begin{prop}
The bilinear form over reals
\begin{center}
(f(x),g(x))$_{ \psi} $ := [(Wf) \textbf{(} $Q\left( {\partial _{\psi} }
\right)$\textbf{)}$S_{\partial _{\psi} }  \;$g(x)]$_{x=0}$ f,g $ \in P$
\end{center}
is a positive definite inner product.
\end{prop}

\begin{proof}
\begin{multline*}
\left( {s_{k} \left( {x} \right)\,,\, s_{n} \left( {x} \right)} \right)
_{\psi} = \textrm{\bf [} Q\left( {\partial _{\psi} }  \right)^{k}S_{\partial
_{\psi}}  \;s_{n} \left( {x} \right) \textrm{\bf ]}_{x=0} = \textrm{\bf [}Q
\left( {\partial _{\psi} }  \right)^{k}p_{n} \left( {x} \right)$\textbf{]}$
_{x=0} = \\ = n^{\underline {k}}_{\psi}  p_{n - k} \left( {0} \right)
= n_{\psi}!  \delta _{nk}
\end{multline*}
where $\left\{ {p_{n} \left( {x} \right)}
\right\}_{n = 0}^{\infty}  $ is $\partial _{\psi}  $-basic polynomial
sequence of the operator $Q\left( {\partial _{\psi} }
\right)$.
\end{proof}
We shall call the scalar product ( , )$_{ \psi} $ : $\aleph
\times \aleph  \to $ \textbf{\textit{R}} - the natural inner product
associated with the sequence $\left\{ {s_{n} \left( {x} \right)} \right\}_{n
= 0}^{\infty}  $ of Sheffer $\psi $-polynomials. Naturally unitary space
(\textit{P} ; ( , )$_{ \psi} $ ) may be completed to the unique Hilbert
space $\aleph = \overline {P} $. If so then the following Theorem is
valid also in $\psi $-extended case of finite operator calculus.

\begin{thm}\label{thfourfive} {\em (Spectral Theorem)}\\
Let $\left\{ {s_{n} \left( {x} \right)} \right\}_{n = 0}^{\infty}  $ be the
sequence of Sheffer $\psi$-polynomials relative to an invertible
$\partial _{\psi}  $-shift invariant operator $S_{\partial _{\psi} }  \;$for
the $\partial _{\psi}  $-delta operator $Q\left( {\partial _{\psi} }
\right)$ with the $\partial _{\psi}  $-basic polynomial sequence $\left\{
{p_{n} \left( {x} \right)} \right\}_{n = 0}^{\infty}  $. Then there exists a
unique essentially self adjoint operator $A_{\psi}  $ : $\aleph  \to $
$\aleph $ given by
\[
A_{\psi}  = \sum\limits_{k \ge 1} {\frac{{u_{k} + \hat {v}_{k} \left( {x}
\right)}}{{\left( {k - 1} \right)_{\psi}  !}}}
Q\left( {\partial _{\psi} }  \right)^{k}
\]
such that the spectrum of $A_{\psi}  $ consists of $\; n = 0, 1, 2, 3,...$
where $A_{\psi}  s_{n} \left( {x} \right) = ns_{n} \left( {x}\right)$.
The quantities $u_{k} $ and $\hat {v}_{k} $ are calculated according to :
\begin{center}
 $u_{k} \; =
\; - \left[ {\left( {logS_{\partial _{\psi} } }  \right)^{`}\hat {x}_{\psi}
^{ - 1}p_{k} \left( {x} \right)} \right]_{x = 0} $ and $\; \hat {v}\left( {x}
\right)_{k}  = \hat {x}_{\psi}  \left[ {\frac{{d}}{{dx}}p_{k}^{}}
\right]\left( {0} \right)$.
\end{center}
\end{thm}

\begin{proof}
>From the second expansion theorem one has
\[
S_{\partial _{\psi} }  ^{ - 1}\;
E^{a}\left( {\partial _{\psi} }  \right) \quad =
\sum\limits_{n \ge 0} {\frac{{s_{n} \left( {a} \right)}}{{n_{\psi}  !}}}
Q\left( {\partial _{\psi} ^{}}  \right)^{n}.
\]
Taking the Pincherle derivative of both sides in the above we get

\begin{center}
($S_{\partial _{\psi} }  ^{ - 1}\;E^{a}\left( {\partial _{\psi} }
\right)$)\textbf{'} =$\sum\limits_{n \ge 1} {\frac{{s_{n} \left( {a}
\right)}}{{n_{\psi}  !}}} nQ\left( {\partial _{\psi} ^{}}  \right)^{n -
1}Q\left( {\partial _{\psi} }  \right)$\textbf{'}.
\end{center}
After multiplication of both sides by ($Q\left( {\partial _{\psi} }
\right)$\textbf{'} )$^{-1}Q\left( {\partial _{\psi} }  \right)$ we arrive
at
\begin{multline}\label{star}
\left( {-S_{\partial _{\psi} }  ^{ - 1}\;S_{\partial _{q}}  \; + \hat {a}_{\psi
}} \right) S_{\partial _{\psi} }  ^{ - 1}\;E^{a}\left( {\partial _{\psi} }
\right) \left( {Q\left( {\partial _{\psi} }  \right) \textrm{\bf '}}
\right)^{-1}Q\left(
{\partial _{\psi} }  \right) \equiv \\
 \equiv  \quad T_{\partial _{\psi} }  \;S_{\partial _{\psi} }  ^{ -
1}\;E^{a}\left( {\partial _{\psi} }  \right) = \sum\limits_{n \ge 1}
{\frac{{s_{n} \left( {a} \right)}}{{n_{\psi}  !}}} nQ\left( {\partial_{\psi}}  \right)^{n}
\end{multline}
where
\begin{multline*}
 T_{\partial _{q}}  \; \equiv  (-S_{\partial _{\psi} }  ^{ -
1}\;S_{\partial _{\psi} }  \;$\textbf{'} $+ \hat {a}_{\psi})(Q\left(
{\partial _{\psi} }  \right)$\textbf{'} )$^{-1}Q\left( {\partial _{\psi}
} \right) = (\hat {a}_{\psi} - ($log$S_{\partial _{\psi} }  \;$\textbf{'})
($Q\left( {\partial _{\psi} }  \right)$\textbf{'} )$^{-1}Q\left( {\partial
_{\psi} }  \right)
\end{multline*}
Using the Rodrigues formula $p_{n+1}(x) = \frac{{\left(
{n + 1} \right)_{\psi} } }{{n + 1}}\hat {x}_{\psi} \left(
{Q\left( {\partial_{\psi} }  \right) \textrm{\bf '}} \right)^{-1}
p_{n}$\textit{(x)}  we define $q_{n} \left( {x} \right)$:
\[
q_{n} \left( {x} \right): \quad =
\quad
\frac{{\left( {n + 1} \right)}}{{\left( {n + 1} \right)_{\psi} } }
\hat {x}_{\psi}  ^{ - 1} p_{n+1}(x)  = \left( {Q\left(
{\partial _{\psi} }  \right) \textrm{\bf'}} \right)^{-1}p_{n}(x)
\]
for $n \ge 0$ whence
\[
Q\left( {\partial _{\psi} }  \right)
q_{n} \left( {x} \right) \quad = \left( {Q\left( {\partial _{\psi} }\right)
\textrm{\bf'}} \right)^{-1} \quad \left( {\partial _{\psi} }
\right)p_{n}(x) = n_{\psi}  q_{n - 1} \left( {x} \right)
\]
for $n>0$.
Consider now $\quad T_{\partial _{\psi} }  = \sum\limits_{n \ge 0}
{\frac{{\hat{b}_{n}} }{{n_{\psi}  !}}} Q\left( {\partial _{\psi}}
\right)^{n} \quad$
where $\quad \hat {b}_{k} = \left. {\left[ {T_{\partial _{\psi} }  p_{k}
\left( {x} \right)} \right]} \right|_{x = 0}$. Easy calculation yields
for $k \in N$ : $\hat {b}_{k} = \left. {\left[ {T_{\partial _{\psi} }  p_{k} \left( {x}
\right)} \right]} \right|_{x = 0} =  k \left(\hat {v}_{k}  + u_{k} \right)$
where
\[
u_{k} = \;
 - \left[
 {\left(
 {\textrm{log}S_{\partial _{\psi} } }
   \right)^{`}\hat {x}_{\psi}
^{ - 1}p_{k} \left( {x} \right)}
 \right]_{x = 0} \textrm{and } \hat {v}_{k} \left(
{a} \right) = \left[
{\hat {a}_{\psi}  \hat {x}_{\psi}  ^{ -
1}p_{k} \left( {x} \right)
} \right]_{x = 0} = \hat {a}_{\psi}  \left[
{\frac{d}{dx}p_{k}}  \right]\left( {0} \right).
\]
In order to find out the operator $A_{\psi}  $ ; $A_{\psi}  s_{n} \left(
{x} \right) = ns_{n} \left( {x} \right)$ it is enough to consider
now an expansion of the operator $T_{\partial _{\psi} }  \;S_{\partial
_{\psi} }  ^{ - 1}\;E^{a}\left( {\partial _{\psi} }  \right)$. For that
purpose note that
\[
T_{\partial _{\psi} }  \;
S_{\partial _{\psi} }  ^{ - 1}\;
f\left( {x + _{\psi}  a} \right)=
\sum\limits_{k \ge 1} {\frac{{\hat {b}_{k}} }{{k_{\psi}  !}}} Q\left(
{\partial _{\psi} ^{}}  \right)^{k}\left[ {S_{\partial _{q}}  ^{ - 1}\;
f\left({x + _{\psi}  a} \right)} \right] \forall f \in P
\]
and recall the formula we have started the proof with doe to which we
have $S_{\partial _{\psi} }  ^{ - 1}\;E^{a}\left( {\partial _{\psi} }
\right)$
$$S_{\partial _{\psi} }  ^{ - 1}\;f\left( {x + _{\psi}  a}
\right) = \sum\limits_{n \ge 0} {\frac{{s_{n} \left( {x} \right)}}{{n_{\psi
} !}}} Q\left( {\partial _{\psi} ^{}}  \right)^{n}f\left( {a} \right).$$
After the obvious inserting we obtain
\begin{center}
 $T_{\partial _{\psi} }  \;
S_{\partial _{\psi} }  ^{ - 1}\;
f\left( {x + _{\psi}  a} \right)=
\sum\limits_{k \ge 1} {\frac{{\hat {b}_{k}} }{{k_{\psi}  !}}} Q\left(
{\partial _{\psi} ^{}}  \right)^{k} \left[ {\sum\limits_{n \ge 0}
{\frac{{s_{n} \left( {a} \right)}}{{n_{\psi}  !}}} Q\left( {\partial_{\psi}}
\right)^{n}f\left( {a} \right) } \right] \quad \forall \; f \in P$
or
\end{center}
\begin{center}
$T_{\partial _{\psi} }  \;
S_{\partial _{\psi} }  ^{ - 1}\;
f\left( {x + _{\psi}  a} \right)=
\sum\limits_{n \ge 0} \left[ {\sum\limits_{k \ge 1} {\frac{{\hat {b}_{k}
}}{{k_{\psi}  !}}s_{n} \left( {x} \right)} Q\left( {\partial _{\psi} ^{}}
\right)^{k}} \right]
\frac{{Q\left( {\partial _{\psi} }  \right)^{n}}}{{n_{\psi}  !}}
f\left( {a} \right) \quad \forall \; f \in P$.
\end{center}
Permuting up there $a$ and $x$ one gets
\[
T_{\partial _{\psi} }  \;
S_{\partial _{\psi} }  ^{ - 1}\;
E^{a}\left( {\partial _{\psi} }  \right)=
\quad
\sum\limits_{n \ge 0} \left[ {\sum\limits_{k \ge 1} {\frac{{\hat {b}_{k}
}}{{k_{\psi}  !}}s_{n} \left( {a} \right)} Q\left( {\partial _{\psi} ^{}}
\right)^{k}} \right]
\frac{{Q\left( {\partial _{\psi} }  \right)^{n}}}{{n_{\psi}  !}}.
\]
Comparing now this with [see: (\ref{star})]
\[
T_{\partial _{\psi} }  \;
S_{\partial _{\psi} }  ^{ - 1}\;
E^{a}\left( {\partial _{\psi} }  \right)=
\sum\limits_{n \ge 1} {\frac{{s_{n} \left( {a} \right)}}{{n_{\psi}  !}}}
nQ\left( {\partial _{\psi}}  \right)^{n}
\]
we conclude that upon chanching again $a$ to $x$ we end up with
\[
\sum\limits_{k \ge 1} {\frac{{\hat {b}_{k}} }{{k_{\psi}  !}}} Q\left(
{\partial _{\psi} ^{}}  \right)^{k}
s_{n} \left( {x} \right) \; = \; ns_{n} \left( {x} \right) \quad
\textrm{for }n \ge 0 \quad \left( {b_{0}} = 0 \right)
\]
where $\hat {b}_{k} = \left. {\left[ {T_{\partial _{\psi} }  p_{k} \left(
{x} \right)} \right]} \right|_{x = 0} = k \left( {\hat {v}_{k} +
u_{k}} \right) \;$ , $k \in N$.
The operator $A_{\psi}  $ : $\aleph  \to  \quad \aleph $ ; $A_{\psi}  s_{n}
\left( {x} \right)$ = \textit{n}$s_{n} \left( {x} \right)$ has then the form
\[
A_{\psi}  = \sum\limits_{k \ge 1} {\frac{{u_{k} + \hat {v}_{k} \left( {x}
\right)}}{{\left( {k - 1} \right)_{\psi}  !}}}
Q\left( {\partial _{\psi} }  \right)^{k} \; .
\]
The sequence $\left\{ {s_{n} \left( {x} \right)} \right\}_{n =
0}^{\infty}  $of Sheffer $\psi $-polynomials spans the Hilbert space $\aleph
$ We therefore state that the operator $A_{\psi}  $ :\textbf{ $\aleph  \to
$} $\aleph $ is the unique essentiall self adjoint unbounded operator
having nonnegative integers as its spectrum with corresponding eigenfuctions
$\left\{ {s_{n} \left( {x} \right)} \right\}_{n = 0}^{\infty}  $.
\end{proof}
To this end let us draw one general conclusion and let us make a link to
Hopf algebra in a form of a remark with indication of references for further
readings.

\begin{conc} {\em
One has experienced the way in which finite operator calculus gains its
natural $\psi $-extended representation. We have seen that the use of Rota`s
calculus notions and ideas of proofs with indispensable changes leads to
this kind of extension.

In this context it should be noted that the notion of generalised Sheffer
$\psi $-polynomials (see Comment \ref{com1} and \cite{29} ) constitutes
the ultimate generalisation of umbral calculus in that sense that this
calculus deals now with \textit{all} \textit{polynomial} sequences ( see:
\cite{29} and \cite{38} ).}
\end{conc}

\begin{rem} {\em (general Hopf algebras` remark)\\
Exactly as in the undeformed case of umbral calculus $^{} $\cite{38} the
generalized translation operator $E^{y}\left( {\partial _{\psi} }
\right) = exp_{\psi}  \{ y\partial _{\psi}  \} $ is \textit{an example} of
comultiplication over the space of polynomials F[x] ( F[x] = \textit{P})
i.e. a map from F[x] to F[x,y] satisfying the conditions of coassociativity
and counicity ; indeed - with notation $E^{y}\left( {\partial _{\psi} }
\right)p\left( {x} \right)=: p\left( {x + _{\psi}  y} \right)$ we have
$p\left( {\left( {x + _{\psi}  y} \right) + _{\psi}  z} \right) = p\left( {x
+ _{\psi}  \left( {y + _{\psi}  z} \right)} \right)$ and generalized
translation operator fixes all constants.

Hence with natural antipode the algebra of polynomials may be seen as a
graded Hopf algebra \cite{28}. Coalgebra maps \cite{38} S are in their turn exactly
umbral substitution maps as defined in \cite{36} - (see there Theorem 5 then and
section 7). The $\psi $-extension of this Theorem 5 from \cite{36} with
corresponding ``$\psi $-proof'' is automatic. In this connection note: there
are other than $E^{y}\left( {\partial _{\psi} }  \right)$ comultiplications
however all these are of the form

 $E^{y}\left( {\Omega}  \right)=
\quad
\sum\limits_{n \ge 0} {\omega _{n} \left( {y} \right)\frac{{\Omega
^{n}}}{{n_{\psi}  !}}} \quad -$ where $\{ \omega _{n} \left( {x} \right)\} _{n
\ge 0} $ is a basic sequence of $\Omega $.}
\end{rem}
A polynomial sequence $\{ s_{n} \left( {x} \right)\} _{n \ge 0} $ is then
said to be a Sheffer \textit{with respect to} $E^{y}\left( {\Omega}
\right)$ comultiplication \textit{and a pair} of $E^{y}\left( {\Omega}
\right)$ invariant operators $Q$ and $S$
(see Proposition \ref{propfourthree}).

At the end this section we give few examples - postponing the systematic
application of the $\psi $-calculus of Rota to the subsequent publication .
Meanwhile let us indicate reference \cite{7}- as the elegant presentation of
umbral calculus using objects of triple meaning.

There \cite{7} a lot of examples and also some $q$-examples are elaborated. The
latter encompass Gegenbauer and Jacobi polynomials hence also Chebyshev
polynomials. See also \cite{10}. For $q$-Abel polynomials see \cite{26a},
\cite{26b,26c}. For
$q$-Hermite polynomials see \cite{9a,9b} while for $q$-Laguerre polynomials
see \cite{39a,39b}. Also quantum $q$-oscillator algebra provides a
natural setting for $q$-Laguerre polynomials and $q$-Hermite polynomials
\cite{40a,40b,40c,41}.

It seems to be interesting and perhaps important question to ask whether
these quantum like models might be $\psi $-extended.

\subsection*{Examples}
\begin{enumerate}
\renewcommand{\labelenumi}{\bf \arabic{enumi}.}

\item Let $Q\left( {\partial _{\psi} }  \right) = \partial _{\psi}  $
then $\partial _{\psi}  $\textbf{'} $= id$. Apply now the Rodrigues
$\psi $-formula in a recurrent way
\begin{multline*}
p_{n}(x)  =  \frac{n_{\psi}}{n} \hat {x}_{\psi}
\left( {\partial_{\psi}\textrm{\bf '}} \right)^{-1} p_{n-1}(x) =
\frac{n_{\psi}}{n} \hat {x}_{\psi} \left( {
\frac{\left( {n - 1}\right)_{\psi}}{\left( {n - 1} \right)} \hat {x}_{\psi}
p_{n-2}(x)   } \right) = ... =\\= \frac{n_{\psi}!}{n!} \left( {
\hat {x}_{\psi}} \right)^{n} \left[ {1} \right] = x^{n}.
\end{multline*}

\item Let $Q\left( {\partial _{q}}  \right)$= $\Delta _{q} : = E\left(
{\partial _{q}}  \right) - id$. Then $\Delta _{q} $\textbf{'} $= \hat
{1}_{q} E^{1}\left( {\partial _{q}}  \right)$ where $\hat {1}_{q} x^{n}$
=$\frac{{\left( {n + 1} \right)}}{{\left( {n + 1} \right)_{q}} }x^{n}$ in
accordance with $\hat {x}_{q} x^{n} = \frac{{\left( {n + 1}
\right)}}{{\left( {n + 1} \right)_{q}} }x^{n + 1}$. Let us observe here that
$\psi$-extended case is covered in this example just by replacement $q \to
\psi $.

Recall now the Rodrigues $q$- formula $p_{n}(x) =
\frac{{n_{q}} }{{n}}\hat {x}_{q} $ ($\Delta _{q} $\textbf{'} )$^{-1}$
\textit{p}$_{n-1}$\textit{(x)} = $\frac{{n_{q}} }{{n}} \quad \hat {x}_{q} E^{
- 1}\left( {\partial _{q}}  \right)\hat {1}_{q} ^{ -
1}$\textit{p}$_{n-1}$\textit{(x)} where $\hat {1}_{q} ^{ - 1}x^{n}$
=$\frac{{n_{q}} }{{n}}x^{n}$ for $n>0$ and $\hat {1}_{q} ^{ - 1}$[1] = 1 in
order to apply this Rodrigues $q$- formula in a recurrent way as follows:

\begin{multline*}
p_{n}(x) = \frac{{n_{q}} }{{n}}\hat {x}_{q} \; E^{-1} \left( {\partial _{q}}
\right)\hat {1}_{q} ^{-1}p_{n-1}(x) = \\
= \; \frac{{n_{q}} }{{n}}\hat {x}_{q}
E^{ - 1}\left( {\partial _{q}}  \right)\hat {1}_{q} ^{ -1} \left(
{\frac{{\left( {n - 1} \right)_{q}} }{{\left( {n - 1} \right)}}\hat {x}_{q}
E^{ - 1}\left( {\partial _{q}}  \right)\hat {1}_{q} ^{ -1} p_{n-2} (x)}
\right) =\\
= \; \frac{{n_{q}} }{{n}} \; \frac{{\left( {n - 1} \right)_{q}
}}{{\left( {n - 1} \right)}}\hat {x}_{q} E^{ - 1}\left( {\partial _{q}}
\right)\hat {1}_{q} ^{ - 1} \left( {E^{ - 1}\left( {\partial _{q}}
\right)\hat {1}_{q} ^{ - 1}p_{n-2}(x)} \right) = ...=\\
= \frac{{n_{q} !}}{{n!}}\{ \left( {\hat {x}_{q} E^{ - 1}\left( {\partial _{q}}
\right)\hat {1}_{q} ^{ - 1}} \right)\circ \left( {\hat {x}_{q} E^{ - 1}\left(
{\partial _{q}}  \right)\hat {1}_{q} ^{ - 1}}\right)\circ  ... \circ
\left ({\hat {x}_{q} E^{ - 1}\left( {\partial _{q}}  \right)\hat {1}_{q} ^{ -
1}}\right)\} \left[{1}\right]  \equiv\\
 \equiv
 \quad \frac{{n_{q} !}}{{n!}}\left({\hat {x}_{q} E^{ - 1}\left(
 {\partial_{q}} \right)\hat {1}_{q} ^{ - 1}}\right)^{n} \left[{1}\right] =
p_{n}(x).
\end{multline*}

Due to $E^{ - 1}\left( {\partial _{q}}  \right) = \sum\limits_{n \ge 0}
{\frac{{\left( { - 1} \right)^{n}}}{{n_{q} !}}} \partial _{q}^{n} $ we have
\[
E^{ - 1}\left( {\partial _{q}}  \right)x^{n} = \sum\limits_{k = 0}^{n}
{} \left( {{\begin{array}{*{20}c}
 {n} \hfill \\
 {k} \hfill \\
\end{array}} } \right)_{q} \left( { - 1} \right)^{k}x^{n - k} \; \equiv
 \left( {\textrm{\em x} -_{q} 1} \right)^{n}.
\]
Hence
$$\left( {\hat {x}_{q} E^{ - 1}\left( {\partial _{q}}
\right)\hat {1}_{q} ^{ - 1}} \right) \left[{1}\right] = x \textrm{ and }$$
$$\left( {\hat {x}_{q} E^{ - 1}\left(
{\partial _{q}}  \right)\hat {1}_{q} ^{ - 1}} \right)\left[{x}\right] =\\
= \; \frac{{2}}{{2_{q}}} x^{2}- x
\mathrel{\mathop{\kern0pt\longrightarrow}\limits_{q \to 1}}
x(x-1) = x^{\underline {2}}$$
and by induction we get the expression:
$$ p_{n}(x) = \frac{{n_{q} !}}{{n!}} \left( {\hat {x}_{q}
E^{ - 1}\left( {\partial _{q}}  \right)\hat {1}_{q} ^{ - 1}} \right)^{n}
\left[{1}\right]$$
which may be given a simpler form
$$p_{n}(x) =
\frac{{n_{q} !}}{{n!}}\hat {x}_{q} \left( {E^{ - 1}\left( {\partial _{q}}
\right)\hat {x}} \right)^{n} \left[{1}\right] \textrm{due to } \hat {1}_{q}
^{ - 1}\hat {x}_{q} = \hat {x}.$$
Of course $p_{n}(x) = \frac{{n_{q} !}}{{n!}}\hat
{x}_{q} \left({E^{ - 1}\left( {\partial _{q}}  \right)\hat {x}}\right)^{n}
$[1]$\mathrel{\mathop{\kern0pt\longrightarrow}\limits_{q \to 1}} x^{\underline
{n}} $ where
$\left\{ {x^{\underline {n}} } \right\}_{n \ge 0} $ is the basic polynomial
sequence of the delta operator $\Delta $ ; $\Delta _{q} $
$\mathrel{\mathop{\kern0pt\longrightarrow}\limits_{q \to 1}} \; \Delta $.\\
Let us observe here that $\psi $-extended case is covered in this example
just by replacement $q \to \psi $.\\
One may give also the explicit expression of $\partial _{q} $-basic
polynomial sequence $\left\{ {p_{n} \left( {x} \right)} \right\}_{n =
0}^{\infty}  $ of the $\partial _{q} $-delta operator $\Delta _{q} $ using
the statement (3) \textit{p}$_{n}$\textit{(x) =} $\frac{{n_{q}} }{{n}}\hat
{x}_{q} $  $S_{\partial _{q}}  \;x^{n-1}$ of the
Theorem~\ref{thfourone}. Namely - we see that
\[
\Delta _{q} = \partial _{q} S_{\partial _{q}}  \; \equiv \partial
_{q} \sum\limits_{k \ge 0} {\frac{{\partial _{q} ^{k}}}{{\left( {k + 1}
\right)_{q} !}}}  = E^{a}\left( {\partial _{q}}  \right)- id \textrm{ i.e.}
S_{\partial _{q}}  \; = \sum\limits_{k \ge 0} {\frac{{\partial _{q}
^{k}}}{{\left( {k + 1} \right)_{q} !}}}
\]
and for $n>0$
\[
p_{n}(x) = \frac{{n_{q}} }{{n}}\hat {x}_{q} \sum\limits_{k \ge 0}
{\frac{{\partial _{q} ^{k}}}{{\left( {k + 1} \right)_{q} !}}}
x^{n-1} = \frac{{n_{q}} }{{n}}\hat {x}_{q} \sum\limits_{k \ge
0}^{n - 1} {\frac{{\left( {n - 1} \right)_{q} ^{\underline {k}}
}}{{\left( {k + 1} \right)_{q}!} }} x^{n - k - 1}
\]
and finally
\begin{center}
\textit{p}$_{n}$\textit{(x)} = $\frac{{n_{q}} }{{n}}\sum\limits_{k
\ge 0}^{n - 1} {\frac{{\left( {n - 1} \right)_{q} ^{\underline
{k}} }}{{\left( {k + 1} \right)_{q}!} }} \frac{{\left( {n - k}
\right)}}{{\left( {n - k} \right)_{q}} }x^{n - k}$ .
\end{center}
Let us observe again that $\psi $-extended case is covered in this example
just by replacement $q$ $ \to \psi $.

\item Let $Q\left( {\partial _{q}}  \right) = \nabla _{q} : = id -
E^{ - 1}\left( {\partial _{q}}  \right)$ . Then $\nabla _{q} $\textbf{'}$=
\hat {1}_{q} E^{ - 1}\left( {\partial _{q}}  \right)$ . Similarly as in
example (2) $p_{n}(x) = \frac{{n_{q} !}}{{n!}} \left( {\hat
{x}_{q} E\left( {\partial _{q}}  \right)\hat {1}_{q} ^{ - 1}} \right)^{n}
\left[{1}\right]$
or due to $\hat {1}_{q} ^{ - 1}\hat {x}_{q} = \hat {x}$ :
$p_{n}(x) = \frac{{n_{q} !}}{{n!}}\hat {x}_{q}
\left({E\left( {\partial _{q}}  \right)\hat {x}}\right)^{n} $[1].

Of course $p_{n}(x) = \frac{{n_{q} !}}{{n!}}\hat
{x}_{q} \left({E\left( {\partial _{q}}  \right)\hat {x}}\right)^{n} $[1]
$\mathrel{\mathop{\kern0pt\longrightarrow}\limits_{q \to 1}}  \quad x^{\overline
{n}} $ where
$\left\{ {x^{\overline {n}} } \right\}_{n \ge 0} $ is the basic polynomial
sequence of the delta operator $\nabla $; $\nabla _{q} $
$\mathrel{\mathop{\kern0pt\longrightarrow}\limits_{q \to 1}} \; \nabla $ .

Let us observe then that $\psi $-extended case is covered in this example
just by replacement $q \to \psi $.

One may give also the explicit expression of $\partial _{q} $-basic
polynomial sequence $\left\{ {p_{n} \left( {x} \right)} \right\}_{n =
0}^{\infty}  $ of the $\partial _{q} $-delta operator $\nabla _{q} $ using
the statement
\[
p_{n}(x) = \frac{{n_{q}} }{{n}}\hat
{x}_{q} S_{\partial _{q}} \; x^{n-1}.
\]
Namely

 $\nabla _{q} \; =
\;
\partial _{q}
S_{\partial _{q}}  \; \equiv \;
\partial _{q} \sum\limits_{k \ge 0} {\frac{{\left( { - 1}
\right)^{k}\partial _{q} ^{k}}}{{\left( {k + 1} \right)_{q} !}}} \; =
\; id - E^{ - 1}\left( {\partial _{q}}  \right)$ i.e. $S_{\partial _{q}}
\;= \; \sum\limits_{k \ge 0} {\frac{{\left( { - 1} \right)^{k}\partial _{q}
^{k}}}{{\left( {k + 1} \right)_{q} !}}} $

and for $n>0$

\textit{p}$_{n}$\textit{(x) =} $\frac{{n_{q}} }{{n}}\hat {x}_{q} $
$\sum\limits_{k \ge 0} {\frac{{\left( { - 1} \right)^{k}\partial
_{q} ^{k}}}{{\left( {k + 1} \right)_{q} !}}} x^{n-1} =
\frac{{n_{q} }}{{n}}\hat {x}_{q} \sum\limits_{k \ge 0}^{n - 1}
{\left( { - 1} \right)^{k}\frac{{\left( {n - 1} \right)_{q}
^{\underline {k}} }}{{\left( {k + 1} \right)_{q}!} }} x^{n - k -
1}$ . Finally

\begin{center}
\textit{p}$_{n}$\textit{(x)=} $\frac{{n_{q}} }{{n}}\sum\limits_{k
\ge 0}^{n - 1} {\left( { - 1} \right)^{k}\frac{{\left( {n - 1}
\right)_{q} ^{\underline {k}} }}{{\left( {k + 1} \right)_{q}!} }}
\frac{{\left( {n - k} \right)}}{{\left( {n - k} \right)_{q}} }x^{n
- k}$.
\end{center}

Naturally $\psi $-extended case is covered in this example just by
replacement $q \to \psi $.

\item Let $A\left( {\partial _{q}}  \right)$= $\partial _{q}
E^{a}\left( {\partial _{q}}  \right)$ be the ``$q$-Abel operator''. Recall
now from the Theorem~\ref{thfourone}. the statement (3) $p_{n}(x) =
\frac{{n_{q}} }{{n}}\hat {x}_{q} S^{-n}\textrm{x}^{n-1}$
for $n>0$ where $\left\{ {p_{n} \left( {x} \right)} \right\}_{n = 0}^{\infty}
$ is the $\partial _{q} $-basic polynomial sequence of the $\partial _{q}
$-delta operator $Q\left( {\partial _{q}}  \right)$ of the form : $Q\left(
{\partial _{q}}  \right) = \partial _{q} S_{\partial _{q}}$. In the case of
$q$-Abel operator we then have for $n>0$

\begin{center}
$A_{n,q}(x) = \frac{{n_{q}} }{{n}}\hat {x}_{q}
E^{ - na}\left( {\partial _{q}}  \right)x^{n-1}$
\end{center}

where we shall call the $\partial _{q} $-basic polynomial sequence $\left\{
{A_{n,q} \left( {x} \right)} \right\}_{n \ge 0} $ of the $q$-Abel operator -
Abel $q$-polynomials sequence. With our convention we may write the formula
for Abel $q$-polynomials in the form which mimics the undeformed case form

\begin{center}
$A_{n,q}(x) = \frac{{n_{q}} }{{n}}\hat {x}_{q} (x -_{q}\, na)^{n-1}$
$\mathrel{\mathop{\kern0pt\longrightarrow}\limits_{q \to 1}} $\textit{
A}$_{n}$\textit{(x)}$ = (x - na)^{n-1}$
\end{center}

where ($x -_{q}$ \textit{na})$^{n-1} \quad  \equiv  \quad \sum\limits_{k =
0}^{n -1} \left( {{\begin{array}{*{20}c}
 {n - 1} \hfill \\
 {k} \hfill \\
\end{array}} } \right)_{q} \left( { - na} \right)^{k}x^{n - k - 1}$ or
explicitly written

\begin{center}
$A_{n,q}(x) = \frac{{n_{q}} }{{n}}\sum\limits_{k =
0}^{n - 1} {} \left( {{\begin{array}{*{20}c}
 {n - 1} \hfill \\
 {k} \hfill \\
\end{array}} } \right)_{q} \left( { - na} \right)^{k}\frac{{\left( {n -
k} \right)}}{{\left( {n - k} \right)_{q}} }x^{n - k}$ .
\end{center}

Again note: $\psi $-extended case is covered in this example just by
replacement $q \to  \psi $.

\item  Let $L\left( {\partial _{q}}  \right) \; = \; - \sum\limits_{k =
0}^{\infty}  {\partial _{q} ^{k + 1}}  \; \equiv \;
\frac{{\partial _{q}}}{{\partial _{q} - 1}} \equiv  -\left[
{\partial _{q} + \partial _{q} ^{2} + \partial _{q} ^{3} +
\partial _{q} ^{4} + \partial _{q} ^{5} +...} \right]$ be the
``$q$-Laguerre operator''. In this case for $n>0$
\begin{multline*}
L_{n,q}(x) = \frac{{n_{q}} }{{n}}\hat {x}_{q}
\left[{\frac{{1}}{{\partial _{q} - 1}}}\right]^{-n}x^{n-1} =
\frac{{n_{q}} }{{n}}\hat {x}_{q} \left( {\partial _{q} - 1}
\right)^{n}x^{n-1} = \\ = \frac{{n_{q}} }{{n}}\hat {x}_{q}
\sum\limits_{k = 1}^{n} \left( { - 1} \right)^{k}\left(
{{\begin{array}{*{20}c}
 {n} \hfill \\
 {k} \hfill \\
\end{array}} } \right)_{q} \partial _{q} ^{n - k}x^{n - 1} =\\
=\frac{{n_{q}} }{{n}}\sum\limits_{k = 1}^{n} {} \left( { - 1}
\right)^{k}\left( {{\begin{array}{*{20}c}
 {n} \hfill \\
 {k} \hfill \\
\end{array}} } \right)_{q} \left( {n - 1} \right)_{q} ^{\underline {n - k}
}\frac{{k}}{{k_{q}} }x^{k} = \frac{{n_{q}} }{{n}}\sum\limits_{k =
1}^{n} {} \left( { - 1} \right)^{k}\frac{{n_{q} !}}{{k_{q}
!}}\frac{{\left( {n - 1} \right)_{q} ^{\underline {n - k}} }}{{\left( {n -
k} \right)_{q} !}}\frac{{k}}{{k_{q}} }x^{k}.
\end{multline*}
So finally

\begin{center}
$L_{n,q}(x) = \frac{{n_{q}} }{{n}}\sum\limits_{k =
1}^{n} {} \left( { - 1} \right)^{k}\frac{{n_{q} !}}{{k_{q} !}}\left(
{{\begin{array}{*{20}c}
 {n - 1} \hfill \\
 {k - 1} \hfill \\
\end{array}} } \right)_{q} \frac{{k}}{{k_{q}} }x^{k}$.
\end{center}

We shall call the $\partial _{q} $-basic polynomial sequence$\left\{
{L_{n,q} \left( {x} \right)} \right\}_{n \ge 0} $of the $q$-Laguerre operator
$L\left( {\partial _{q}}  \right)$ - Laguerre $q$-polynomials sequence- and
note again that $\psi $-extended case is covered in this example just by
replacement $q \to  \psi $.
\end{enumerate}
{\bf Final Remark} The way is now opened to find out $\psi $-analogues of
typical identities derived from explicit form of the above $\partial _{\psi
} $-basic polynomial sequences due to the Corollary~\ref{corfourtwo}. as in
\cite{42} (p.147 - for example Lah numbers) . Also $q$-analogues or $\psi
$-analogues of Hermite and Pollaczek polynomials apart from Laguerre
polynomials are not difficult to be find out so as to give rise to orthogonal
over an interval of real line sequences $\left\{ {s_{n} \left( {x} \right)}
\right\}_{n = 0}^{\infty}  $ of Sheffer $\psi$-polynomials with appropriately
chosen $\partial _{\psi}  $-shift invariant invertible $S_{\partial _{\psi}
} \;$operators.

As stated earlier - in this section we only give few examples - postponing
the systematic application of the $\psi $-calculus of Rota to the subsequent
publication.

Let us then recall again that quantum $q$-oscillator algebra provides a
natural setting for Laguerre $q$-polynomials and Hermite $q$-polynomials
\cite{40a,40b,40c,41}.


\section{No $\psi$-analogue extension of quantum $q$-plane formulation?}

The idea to use ``$q$-commuting variables'' goes back at least to
Cigler (1979) \cite{9a,39a,39b} ( see formula (7),(11) in \cite{9a} ) and
also to Kirchenhofer - see \cite{9b} for further systematic development.
In \cite{9b} Kirchenhofer defined the polynomial sequence
$\left\{ {p_{n}}\right\}_{o}^{\infty}  $ of $q$-binomial type by
$$
p_{n} \left( {A + B} \right) \equiv \sum\limits_{k \ge 0} {\left(
{{\begin{array}{*{20}c}
 {n} \hfill \\
 {k} \hfill \\
\end{array}} } \right)} _{q} p_{k} \left( {A} \right)p_{n - k} \left( {B}
\right)\textrm{ where } [B,A]_{q} \equiv BA-qAB = 0.
$$
$A$ and $B$ might be interpreted here as co-ordinates on quantum $q$-plane
(see \cite{28} Chapter 4). For example $A = \hat {x}$ and $B = y \hat {Q}$
where $\hat {Q}\varphi \left( {x} \right) = \varphi \left( {qx} \right)$ . If
so then the following identification takes place
\[
p_{n} \left( {x + _{q} y} \right) \equiv
E^{y}\left( {\partial _{q}}  \right)
p_{n} \left( {x} \right)=
\sum\limits_{k \ge 0} {\left( {{\begin{array}{*{20}c}
 {n} \hfill \\
 {k} \hfill \\
\end{array}} } \right)} _{q} p_{k} \left( {x} \right)p_{n - k} \left( {y}
\right)=
p_{n} \left( {\hat {x} + y\hat {Q}} \right) \textbf{1}.
\]
$q$-Sheffer polynomials $\left\{ {s_{n} \left( {x} \right)} \right\}_{n =
0}^{\infty}  $ are defined correspondingly by (see: 2.1.1. in \cite{9b})
$s_{n} \left( {A + B} \right) \equiv \sum\limits_{k \ge 0} {\left(
{{\begin{array}{*{20}c}
 {n} \hfill \\
 {k} \hfill \\
\end{array}} } \right)} _{q} s_{k} \left( {A} \right)p_{n - k} \left( {B}
\right)$ where [B,A]$_{q} \equiv BA-qAB = 0$ and the polynomial sequence
$\left\{ {p_{n} \left( {x}
\right)} \right\}_{n = 0}^{\infty}  $ is of $q$-binomial type. For example $A
= \hat {x}$ and $B = y\hat {Q}$ where $\hat {Q}\varphi \left( {x} \right)
= \varphi \left( {qx} \right)$. Then the following identification takes
place:
\[
s_{n} \left( {x + _{q} y} \right) \equiv
E^{y}\left( {\partial _{q}}  \right)
s_{n} \left( {x} \right)=
\sum\limits_{k \ge 0} {\left( {{\begin{array}{*{20}c}
 {n} \hfill \\
 {k} \hfill \\
\end{array}} } \right)} _{q} s_{k} \left( {x} \right)p_{n - k} \left( {y}
\right)=
s_{n} \left( {\hat {x} + y\hat {Q}} \right) \textbf{1}
\]

This means that one may formulate q-extended finite operator
calculus with help of the ``quantum q-plane'' q-commuting variables $A, B:\;
AB - qBA  \equiv [A,B]_{q}= 0$.

The above identifications of polynomial sequences $\left\{ {p_{n}}
\right\}_{o}^{\infty}  $ of $q$-binomial type and Sheffer $q$-polynomials
$\left\{ {s_{n} \left( {x} \right)} \right\}_{n = 0}^{\infty} $ fail to be
extended to the more general $\psi $-case. This means that we can not
formulate that way the $\psi $-extended finite operator calculus with help
of the ``quantum $\psi $-plane'' $\hat {q}_{\psi ,Q} $-commuting variables
$A, B: AB - \hat {q}_{\psi ,Q}BA \equiv [A,B]_{\hat {q}_{\psi ,Q}} = 0$.
We shall explain this in the sequel. For that to do we introduce the following
notions - important on their own.

\begin{defn} \label{defnfiveone}
Let $\left\{ {p_{n}}  \right\}_{n \ge 0} $ be the $\partial _{q} $-basic
polynomial sequence of the $\partial _{q} $-delta operator$Q\left( {\partial
_{q}}  \right)$. A linear map $x_{Q\left( {\partial _{q}}  \right)}
$: $P \to P$ , $x_{Q\left( {\partial _{q}}  \right)} p_{n} = p_{n +
1} ,\quad n \ge 0$ is called the operator dual to $Q\left( {\partial _{q}}
\right)$.
\end{defn}
For $Q = id$ we have: $x_{Q\left( {\partial _{q}}  \right)}
 \equiv  \quad x_{\partial _{q}}   \quad  \equiv  \quad \hat {x}$.\\
\textbf{Comment:} The addjective dual in the above sense corresponds to the
addjective adjoint in $q$-umbral calculus language of linear functionals'
umbral algebra (see : Proposition 1.1.21 in \cite{9b} ). In this connection
we note that in the formulation of undeformed
umbral calculus \cite{10} in terms of algebra \textit{P}$^{*}$ of linear
functionals on \textit{P} - the map $x_{Q\left( {\partial _{\psi} }
\right)} $ is called an ``umbral shift operator'' and it is \textit{adjoint}
to a derivation of \textit{P}$^{*}$ ; see: Theorem 5 in \cite{10}.
See also 1.1.16 in \cite{9b}.

\begin{defn}
Let $\left\{ {p_{n}}  \right\}_{n \ge 0} $ be the $\partial _{\psi}$-basic
polynomial sequence of the $\partial _{\psi}  $-delta
operator $Q\left( {\partial _{\psi} }  \right) = Q$.
Then the $\hat {q}_{\psi ,Q} $-operator is a liner map\\ $\hat {q}_{\psi
,Q} :P \to P; \quad \hat {q}_{\psi ,Q} p_{n} = \frac{{\left( {n + 1}
\right)_{\psi}  - 1}}{{n_{\psi} } }p_{n},\quad n \ge 0$.
\end{defn}

We call the $\hat {q}_{\psi ,Q} $ operator the $\hat {q}_{\psi ,Q}
$\textit{-mutator operator}.\\
\textbf{Example:} For $Q =id\;$ $Q\left( {\partial _{\psi} }
\right) = \partial _{\psi}  $ the natural notation is $\hat {q}_{\psi ,id}
 \equiv \hat {q}_{\psi}  $. For $Q =id$ and $\psi _{n} \left( {q}
\right) = \frac{{1}}{{R\left( {q^{n}} \right)!}}$ and $R\left( {x} \right)
= \frac{{1 - x}}{{1 - q}} \quad \hat {q}_{\psi ,Q}  \equiv \hat {q}_{R,id}
 \equiv \hat {q}_{R} \equiv \hat {q}_{q,id}  \equiv \hat
{q}_{q}  \equiv \hat {q}$ and $\hat {q}_{\psi ,Q} x^{n}$ =
$q^{n}x^{n}$.

\begin{defn}
Let $A$ and $B$ be linear operators acting on $P$; $A:P \to P$, $B: P
\to P$. Then $AB - \hat {q}_{\psi ,Q}BA \equiv  [A,B]_{\hat
{q}_{\psi ,Q}}  $ is called $\hat {q}_{\psi ,Q} $-mutator of $A$ and $B$
operators.
\end{defn}

\begin{obs} {\em
$Q\left( {\partial _{\psi} }  \right) x_{Q\left(
{\partial _{\psi} }  \right)} - \hat {q}_{\psi ,Q} x_{Q\left(
{\partial _{\psi} }  \right)} Q\left( {\partial _{\psi} }  \right)
\equiv $ [$Q\left( {\partial _{\psi} }  \right)$,$\hat {x}_{Q\left(
{\partial _{\psi} }  \right)} $]$_{\hat {q}_{\psi ,Q}} = id$.\\
This is easily verified in the $\partial _{\psi}  $-basic $\left\{ {p_{n}}
\right\}_{n \ge 0} $ of the $\partial _{\psi}  $-delta operator $Q\left(
{\partial _{\psi} }  \right)$.

Equipped with pair of operators ( $Q\left( {\partial _{\psi} }  \right)$ ,
$x_{Q\left( {\partial _{\psi} }  \right)} $ ) and $\hat {q}_{\psi ,Q}
$-mutator we have at our disposal all possible representants of ``canonical
pairs'' of differential operators on the \textit{P} algebra. The meaning of the
adjective: ``canonical'' includes also the content of the
Remark~\ref{remfivetwo}.}
\end{obs}
For important historical reasons at first here is the Remark~\ref{remfiveone}.

\begin{rem} \label{remfiveone} {\em
The $\psi $-derivative is a particular example of a linear operator that
reduces by one the degree of any polynomial. In 1901 it was proved \cite{43}
that every linear operator $T$ mapping $P$ into $P$ may be represented as
infinite series in operators $\hat {x}$ and $D$.
In 1986 the authors of \cite{44} supplied us with the explicit expression for
such series in most general case of polynomials in one variable. We quote here
the Proposition 1 from \cite{44} one has:}
\end{rem}

{\em  Let $L$ be a linear operator that reduces by one each
polynomial. Let $\{ q_{n} \left( {\hat {x}} \right)\} _{n \ge 0} $
be an arbitrary sequence of polynomials in the operator $\hat {x}$.
Then $T = \sum\limits_{n \ge 0} {q_{n} \left( {\hat {x}} \right)}$
$L ^{n}$ defines a linear operator that maps
polynomials into polynomials. Conversely, if T is linear operator that
maps polynomials into polynomials then there exists a unique expansion of
the form
 $$T = \sum\limits_{n \ge 0} {q_{n} \left( {\hat {x}} \right)} L ^{n}.$$}
In 1996 this was extended to algebra of many variables polynomials \cite{45}.

\begin{rem}\label{remfivetwo} {\em
The importance of the pair of dual operators : $Q\left( {\partial _{\psi} }
\right)$ and $x_{Q\left( {\partial _{\psi} }  \right)} $ is reflected
by the facts:\\
a) $Q\left( {\partial _{\psi} }  \right) x_{Q\left( {\partial _{\psi
}}  \right)} $- $\hat {q}_{\psi ,Q} x_{Q\left( {\partial _{\psi} }
\right)} Q\left( {\partial _{\psi} }  \right) \equiv  [Q\left(
{\partial _{\psi} }  \right)$, $x_{Q\left( {\partial _{\psi} }
\right)} $]$_{\hat {q}_{R,Q}}  = id$.\\
b)\textit{ Let} $\{ q_{n} \left( {x_{Q\left( {\partial _{\psi} }
\right)}}  \right)\} _{n \ge 0} $\textit{ be an arbitrary sequence of
polynomials in the operator} $\hat {x}_{Q\left( {\partial _{\psi} }
\right)} $.\textit{ Then} $T = \sum\limits_{n \ge 0} {q_{n} \left( {
x_{Q\left( {\partial _{\psi} }  \right)}}  \right)} Q\left( {\partial
_{\psi} }  \right)^{n}$ \textit{defines a linear operator that maps
polynomials into polynomials. Conversely , if T is linear operator that
maps polynomials into polynomials then there exists a unique expansion of
the form}

 $T = \sum\limits_{n \ge 0} {q_{n} \left( {x_{Q\left( {\partial _{\psi
}}  \right)}}  \right)} Q\left( {\partial _{\psi} }  \right)^{n}$.}
\end{rem}

One may now be tempted to formulate the basic notions of $\psi $-extended
finite operator calculus with help of the ``quantum $\psi $-plane'' $\hat
{q}_{\psi ,Q} $-commuting variables $A, B: [A,B]_{\hat {q}_{\psi ,Q}} = 0$
exactly in the same way as in \cite{9a,9b}. For that to try consider
appropriate generalization of $A = \hat {x}$ and $B = y\hat {Q}$ where this
time the action of $\hat {Q}$ on $\left\{ {x^{n}} \right\}_{0}^{\infty}  $
is to be found from the condition $AB - \hat {q}_{\psi}BA \equiv
[A,B]_{\hat {q}_{\psi} } = 0$. Acting with $[A,B]_{\hat {q}_{\psi} }  $
on $\left\{ {x^{n}} \right\}_{0}^{\infty}  $ due to
 $\hat {q}_{\psi}  x^{n} = \frac{{\left( {n + 1} \right)_{\psi}  -
1}}{{n_{\psi} } }x^{n},\quad n \ge 0$ one easily sees that now $\hat
{Q}x^{n}$ = $b_{n} x^{n}$ where $b_{0} = 0$ and $b_{n}
= \prod\limits_{k = 1}^{n} {\frac{{\left( {k + 1} \right)_{\psi} -
1}}{{k_{\psi} } }} $ for $n>0$ is the solution of the difference equation:
\[
b_{n} - b_{n - 1} \frac{{\left( {n + 1} \right)_{\psi}  - 1}}{{n_{\psi} }}
= 0\quad ,\quad n > 0.
\]
With all above taken into account one immediately verifies that for our $A$
and $B$ $\hat {q}_{\psi}  $-commuting variables
\[
\left( {A + B} \right)^{n} \ne \sum\limits_{k \ge 0} {\left(
{{\begin{array}{*{20}c}
 {n} \hfill \\
 {k} \hfill \\
\end{array}} } \right)} _{\psi}  A^{k}B^{n - k}
\]
unless $\psi _{n} \left( {q} \right) = \frac{{1}}{{R\left( {q^{n}}
\right)!}}$ ; $R\left( {x} \right) = \frac{{1 - x}}{{1 - q}}$ hence $\hat
{q}_{\psi ,Q}  \equiv \hat {q}_{R,id}  \equiv \hat {q}_{R}
\equiv  \hat {q}_{q,id}  \equiv \hat {q}_{q}  \equiv \hat {q}$ and
$\hat {q}_{\psi ,Q} x^{n}$ = $q^{n}x^{n}$. Therefore in conclusion one affirms
that the case of $q$-deformed finite operator calculus is fairly enough
distinguished by the Kirchenhofer approach and the quantum plane notion.

\begin{rem}{\em It is therefore not a supprise that the $q$-deformations are
now a daily bread for $q$-theoretician physicists. We have mentioned also an
``intermediate'' $R$-case because it is of primary importance for
advanced theory of coherent states \cite{18} and it seams to be so in
connection with extended binomial theorem recently proved in \cite{18}.}
\end{rem}
Our observation on this occasion is that unless the condition: $\exists $
$\varphi $ ; $R(q^{n + 1}) - \varphi \left( {q}
\right)R(q^{n}) = 1$, where $R(1)=1$ is satisfied neither
the method from \cite{46} nor the method from \cite{47} used in proof of $q$-binomial
theorem is successful. The extension of iterative method works fine
for general $R$-case \cite{18}.

As for the eventual $\psi $-binomial theorem however - the present author
feels helpless with all the methods mentioned.

{\bf Remark on references:}
Of course references are not complete. Nevertheless we would like to indicate
at the end some of them as the source of further references.

Among them the papers \cite{48} by Andrews (1971) and \cite{49} by Goldman and
Rota (1970) are to be quoted as the ones in which $q$-{\underline {umbral}}
calculus is being started and applied at first - to our knowledge.

The very recent source of sources for references on umbral calculus and its
generalisations developed by Loeb and others is ``\textit{The World of
Generating Functions and Umbral Calculus''} (1999) \cite{38}. See also paper
\cite{50} by Loeb and Rota.

\textbf{Acknowledgements}

We highly acknowledged here inspiration by lectures of Professor Viskov and
discussions with him at Bia{\l}ystok University in December 1999.

The author is also very much indebted to the Referee and to the Editor whose
indications allowed preparing the paper in more desirable form.

I also thank warmly Ewa Gr\c{a}dzka for effective cooperation.

\end{document}